%% file: BurgessVinogradovHighDim_V7_arXiv_Jan2016.tex
\documentclass[oneside,11pt]{amsart}
\usepackage{amssymb,latexsym,amsmath,amsthm}
\usepackage{fancyhdr}
\usepackage{color}
\usepackage{mathrsfs}

\usepackage[margin=1.25in]{geometry}

\input{format}

\numberwithin{equation}{section}

	\newcommand{\xtra}[1]{}

\newcommand{\taubf}{\boldsymbol{\tau}}
\newcommand{\xibf}{\boldsymbol{\xi}}
\newcommand{\taumax}{\tau_{\mathrm{max}}}

\newcommand{\pbf}{\mathbf{p}}
\newcommand{\qbf}{\mathbf{q}}

\newcommand{\Cbf}{\mathbf{C}}
\newcommand{\Fbf}{\mathbf{F}}
\newcommand{\Pbf}{\mathbf{P}}

\newcommand{\Fscr}{\mathscr{F}}
\newcommand{\Gscr}{\mathscr{G}}

\begin{document}

\title[Burgess bounds for multi-dimensional sums]{Burgess bounds for multi-dimensional short mixed character sums} 

\author{L. B. Pierce}
\address{Department of Mathematics, Duke University, 120 Science Drive, Durham NC 27708}
\email{pierce@math.duke.edu}



  
  \begin{abstract}
  This paper proves Burgess bounds for short mixed character sums in multi-dimensional settings. The mixed character sums we consider involve both an exponential evaluated at a real-valued multivariate polynomial $f$, and a product of multiplicative Dirichlet characters.  We combine a multi-dimensional Burgess method with  recent results on multi-dimensional Vinogradov Mean Value Theorems for translation-dilation invariant systems in order to prove character sum bounds in $k \geq 1$ dimensions that recapture the Burgess bound in dimension 1. Moreover, we show that by embedding any given polynomial $f$ into an advantageously chosen  translation-dilation invariant system constructed in terms of $f$, we may in many cases significantly improve the bound for the associated character sum, due to a novel phenomenon that occurs only in dimensions $k \geq 2$. 
    \end{abstract}
  
\maketitle

\section{Introduction}
Let $\chi(n)$ be a non-principal multiplicative Dirichlet character to a modulus $q$, and consider the character sum 
\beq\label{SNH}
 S(N,H) = \sum_{N < n \leq N+H} \chi(n).
 \eeq
The P\'{o}lya-Vinogradov inequality states that
\[S(N,H) \ll q^{1/2} \log q,\] which is nontrivial only if the length $H$ of the character sum is longer than $q^{1/2}\log q$. Burgess famously improved on this in a series of papers \cite{Bur57} \cite{Bur62B} \cite{Bur63A} \cite{Bur86}, proving (among more general results) that for
$\chi$ a non-principal multiplicative character to a prime modulus $q$,
\beq\label{Burgess}
 S(N,H) \ll H^{1-\frac{1}{r}}q^{\frac{r+1}{4r^2}}\log q,
 \eeq
for any integer $r\geq 1$, uniformly in $N$.
This provides a nontrivial estimate for $S(N,H)$ as soon as $H >
q^{1/4+\ep}$; more precisely if $H=q^{1/4+\kappa}$, then the Burgess
bound is of size $Hq^{-\del}$ with 
\beq\label{Burgess_del}
\del \approx \kappa^2.
\eeq
The Burgess bound found immediate applications in an upper bound for the least quadratic non-residue modulo a prime and a celebrated sub-convexity estimate for Dirichlet $L$-functions, and has since been used in a wide range of problems in analytic number theory. Burgess's original strategy has also been refined and simplified (for very recent examples see \cite{GaMo10} \cite{HB12}) and adapted to other problems (for example \cite{HB82} \cite{Pie05}), but its main utility currently remains limited to a few types of short character sums.  It would be highly desirable to generalize the Burgess method further to a wide range of character sums involving additive and multiplicative characters, polynomial arguments, and multiple dimensions.

 In the present work we develop  Burgess bounds for multi-dimensional short mixed character sums of the following form.
For each $i=1,\ldots, k$, let $\chi_i$ be a non-principal multiplicative character modulo a prime $q_i$. Let $f$  be a real-valued polynomial of total degree $d$ in $k$ variables and set
\[ S_k(f;\Nbf,\Hbf) = \sum_{\bstack{\xbf \in \Z^k}{\xbf \in (\Nbf, \Nbf+\Hbf]}} e(f(\xbf)) \chi_1(x_1)\cdots \chi_k(x_k)\]
for any $k$-tuple $\Nbf = (N_1,\ldots, N_k)$ of real numbers and $k$-tuple $ \Hbf = (H_1,\ldots, H_k)$ of positive real numbers, where
\[ (\Nbf, \Nbf + \Hbf ] = (N_1, N_1+H_1] \times (N_2, N_2 + H_2]\times \cdots \times (N_k, N_k + H_k]\]
denotes the corresponding box in $\R^k$, with volume $\|\Hbf \| :=H_1 \cdots H_k$. Note that we do not assume the primes $q_i$ are distinct, and in particular an interesting special case arises when all the $q_i$ are equal to a fixed prime $q$. To avoid vacuous cases we always assume $f$ has positive degree with respect to each of the $k$ variables, and that $H_i \geq 1$ for $i=1,\ldots, k$.

We note the trivial bound
\beq\label{S_triv}
 S_k(f;\Nbf,\Hbf)	\ll  \|\Hbf\|.
	\eeq
Nontrivial upper bounds for $S_k(f;\Nbf,\Hbf)$, particularly when $H_i$ is ``short'' relative to $q_i$, are expected to have a variety of applications, for example to counting integral points on certain hypersurfaces, such as multi-dimensional generalizations of the Markoff-Hurwitz and Dwork hypersurfaces (see related work \cite{Shp14a}, \cite{ChSh14a}).

We will prove bounds that are nontrivial when $H_i \gg q_i^{1/4 + \ep}$ by developing a multi-dimensional version of the Burgess method that allows us to apply recent results of Parsell, Prendiville and Wooley \cite{PPW13} on multi-dimensional Vinogradov Mean Value Theorems. The basic framework of this approach is inspired by \cite{HBP14a}, which treats the one-dimensional case, but a new phenomenon arises in dimensions $k \geq 2$. To make this phenomenon clear, we focus now on two specific results which we may frame in very concrete terms. (Both are immediate corollaries of our most general result, Theorem \ref{thm_gen_PPW}, which is stated in terms of translation-dilation invariant systems; see Section \ref{sec_main_thm}.)

The key strategy of our multi-dimensional  Burgess method will transform the original sum $S_k(f;\Nbf,\Hbf)$ into a collection of many shorter sums $S_k(\tilde{f};\tilde{\Nbf},\tilde{\Hbf})$ with other polynomials $\tilde{f}$ and tuples $\tilde{\Nbf},\tilde{\Hbf}$. The transformations $\tilde{f}$ of $f$ will live inside a certain family, which we may choose to construct in various ways.
If we embed $f$ into the family of all polynomials in $k$ variables of degree at most $d$, we obtain a direct generalization of the work of \cite{HBP14a} to $k$ dimensions (Theorem \ref{thm_box_Vin}). 
But a more sophisticated embedding of $f$ into a potentially much smaller family of polynomials allows us to obtain a sharper result (Theorem \ref{thm_gen_PPW_cor}). We now describe these two results.

\subsection{Generic embedding}
We suppose we are given a  fixed real-valued polynomial $f$ of total degree $d$ in $k$ variables, and a corresponding sum $S_k(f;\Nbf,\Hbf)$.
For $\xbf \in \Z^k$ we will use multi-index notation, so that for a tuple $\be = (\be_1, \cdots, \be_k) \in \Z_{\geq0}^k$ we have $\xbf^\be = x_1^{\be_1}\cdots x_k^{\be_k}$. We let $|\be| = \be_1 + \cdots + \be_k$ denote the total degree of the monomial $\xbf^\be$. We consider the system of Diophantine equations given by 
\beq\label{Dioph_system}
\xbf_{1}^\be + \cdots + \xbf_{r}^\be=  \xbf_{r+1}^\be + \cdots + \xbf_{2r}^\be, \quad \text{for all $1 \leq |\be| \leq d$}
\eeq
where each $\xbf_{j} \in \Z^k$.
We let $R_{d,k}$ denote the number of equations in this system  and $M_{d,k}$ denote the sum of the total degrees appearing in the system; we recall that 
\beq\label{RE_dfn_0}
 R_{d,k} =  \binom{k+d}{k} -1, \qquad 
 M_{d,k} =d\binom{k+d}{k}\frac{k}{k+1}.
 \eeq
 
 Let  $J_{r,d,k}(X)$ denote the number of solutions to the system (\ref{Dioph_system}) with $1 \leq x_{j,i} \leq X$ for all $1 \leq j \leq 2r$, $1 \leq i \leq k$.
The main conjecture in the setting of multi-dimensional Vinogradov Mean Value Theorems is that 
for all $r$ sufficiently large with respect to $d$ and $k$,
\beq\label{J_bound}
 J_{r,d,k}(X) \ll X^{2rk - M_{d,k} + \ep}.
 \eeq
In the case $d=1$, (\ref{J_bound}) holds trivially for all $k,r \geq 1$. The case of $d \geq 2$ is highly nontrivial. Nevertheless,
recently Parsell, Prendiville, and Wooley \cite{PPW13} have proved this for a nearly optimal range of $r$:
\begin{letterthm}[Theorem 1.1 of \cite{PPW13}]\label{thm_PPW}
For $k \geq 1$, $d \geq 2$, if $r \geq R_{d,k} (d+1)$,
the main conjecture (\ref{J_bound}) holds for every $\ep>0$.
\end{letterthm}

We combine this with the Burgess method to prove our first result; 
here we use  the conventions that $\qbf = (q_1,\ldots, q_k)$, $\| \qbf \| = q_1 \cdots q_k$,  $q_{\max} = \max \{q_i\}$, $q_{\min} = \min \{q_i\}$.
\begin{thm}\label{thm_box_Vin}
Let $M=M_{d,k}$. For any  $r > R_{d,k}(d+1)$, if $q_i^{\frac{1}{2(r-M)}} <H_i < q_i^{\frac{1}{2} + \frac{1}{4(r-M)}}$ for each $i=1,\ldots, k$, then 
\[
S_k(f;\Nbf,\Hbf)
	\\
	\ll  \| \Hbf\|^{1-\frac{1}{r}} \|\qbf\|^{\frac{-r-M+1}{4r(r-M)} + \frac{M}{4kr(r-M)}+\ep} q_{\max}^{\frac{2rk}{4r(r-M)}} q_{\min}^{-\frac{M}{4r(r-M)}} ,
\]
uniformly in $\Nbf$, with implied constant dependent on $r,d,k,\ep$ and independent of the coefficients of $f$.
\end{thm}
\xtra{We check that $R (d+1) \geq M$ iff
\[ R(d+1) \geq d(R+1)(\frac{k}{k+1}).\]
Since $k/(k+1)<1$, it is sufficient just to check that 
\[ R(d+1) \geq d(R+1).\]
This holds as long as  $R \geq d$. This is stating that 
\[ \binom{k+d}{k}-1 \geq d \]
or equivalently 
\[ (k+d)! \geq k! (d+1)! \]
which can be re-written as 
\[ (d+k)(d-1+k) \cdots (1+k) \geq (d+1) (d-1+1) \cdots (1+1),\]
which is true if $k \geq 1$.
}
It is illustrative to record the result when all the moduli $q_i$ are equal:
\begin{cor}\label{thm_box_Vin_corq}
Under the hypotheses of Theorem \ref{thm_box_Vin}, if in addition $q_1 = \cdots =q_k = q$ for a fixed prime $q$, 
\[
S_k(f;\Nbf,\Hbf)	\ll  \| \Hbf\|^{1-\frac{1}{r}} q^{\frac{k(r+1-M)}{4r(r-M)}+\ep}.
\]
	\end{cor}
In general, for any dimension $k$, we may check the strength of Corollary \ref{thm_box_Vin_corq} as follows: it is nontrivial if each $H_i = q^{1/4 + \kappa_i}$ for $\kappa_i >0$, in which case choosing $r$ optimally shows that Corollary \ref{thm_box_Vin_corq} provides a bound of size 
$ \| \Hbf \| q^{-\del},$
where 
\[
\del \approx  \frac{\left(\sum_{i=1}^k \kappa_i \right)^2}{k}.
\]
Thus the improvement in the bound $\|\Hbf\|q^{-\del}$ over the trivial bound  is  independent of the degree $d$ of the polynomial $f$, and recovers the Burgess result (\ref{Burgess_del}) when $k=1$. Moreover, Theorem \ref{thm_box_Vin}  recovers Theorem 1.3 of \cite{HBP14a} in dimension $k=1$.

Thus this is a natural multi-dimensional generalized Burgess bound for $S_k(f; \Nbf,\Hbf)$. However, in dimensions $k \geq 2$ another effect can come into play, which we now describe.

\subsection{Minimal embedding}
 Given a polynomial $f$ of total degree $d$ in $k$ variables, we can write it in terms of its coefficients $f_\be$ as 
\[ f (\xbf )  =  \sum_{\be \in \Lambda(f)} f_\be \xbf^\be,\]
where $\Lambda (f)$ is the set of nonzero multi-indices corresponding to monomials in $f$ with non-zero coefficients.  (Since the size of $|S_k(f;\Nbf,\Hbf)|$ is unaffected by any constant term in $f$, we may assume that $f$ has no constant term.) 
Next, we construct a set comprised of all the distinct non-constant monomials (rescaled to be monic) that appear in partial derivatives of $f$ of any order; we will call this set $\Fbf(f)$. To construct $\Fbf(f)$ explicitly, define the ordering $\al \leq \be$ for multi-indices $\al,\be \in \Z^k_{\geq 0}$ to mean that $\al_i \leq \be_i$ for each $1 \leq i \leq k$. Then we see that 
\beq\label{Ff_dfn}
 \Fbf(f) = \{ \xbf^\al : \al \neq (0,\ldots,0), \al \leq \be \text{ for some $\be \in \Lambda(f)$}\}.
 \eeq
Note that we assume $\Fbf(f)$ contains only distinct elements.

Clearly, if we define 
\beq\label{F_dfn}
 \Fbf_{d,k} = \{ \xbf^\al : 1 \leq |\al| \leq d\},
 \eeq
then  for any polynomial $f$ of degree $d$,
\beq\label{inclusion}
 \Fbf (f) \subseteq  \Fbf_{d,k}.
 \eeq
In fact, typically $\Fbf(f)$ will be smaller than $\Fbf_{d,k}$.
We define $R(f)$ to be the number of elements in $\Fbf(f)$, and $M(f)$ to be the sum of the total degrees of the elements in $\Fbf(f)$. (We thus see that $R(f), M(f)$ are analogous to $R_{d,k}, M_{d,k}$.)
Finally, we define the multi-index $\ga(f) \in \Z^k_{\geq 0}$ to be the sum of all the multi-indices occurring in $\Fbf(f)$.

{\bf Example A.} For a simple example, if $f$ is itself a monomial, say
\[ f(\xbf) = \xbf^\del,\]
for a fixed multi-index $\del = (d_1,\ldots, d_k)$, then in this case we would have
\beq\label{Ffset}
\Fbf(f) = \{ \xbf^\al : 0 \leq \al_i \leq d_i \text{ for } 1 \leq i \leq k, \al \neq (0,\ldots, 0)\} .
\eeq
Upon defining
\[\Dcal =  \prod_{1 \leq i \leq k} (d_i+1) ,\] 
a simple calculation shows that in this case the set $\Fbf(f)$ has cardinality
\beq\label{Rf1}
 R(f) := \# \Fbf(f)= \left(\sum_{0 \leq \al_1 \leq d_1} \cdots \sum_{0 \leq \al_k \leq d_k} 1\right) -1 = \Dcal -1,
 \eeq
and the sum of the total degrees of the multi-indices in the set $\Fbf(f)$ is
\beq\label{Mf1}
M(f) := \sum_{0 \leq \al_1 \leq d_1} \cdots \sum_{0 \leq \al_k \leq d_k} (\al_1 + \cdots + \al_k)
	= \frac{1}{2} \Dcal   (d_1 + \cdots + d_k) .
	\eeq
\xtra{
\begin{eqnarray*}
M(f)  &=&  \sum_{0 \leq \al_1 \leq d_1} \cdots \sum_{0 \leq \al_k \leq d_k} (\al_1 + \cdots + \al_k) \\
	& = & \sum_{i =1}^k \left( \sum_{0 \leq j \leq d_i} j \right) \left( \prod_{l \neq i}(d_l+1) \right) \\
	& = & \sum_{i =1}^k \frac{d_i(d_i+1)}{2} \left( \prod_{l \neq i}(d_l+1) \right) \\
	& = & \frac{1}{2} (d_1 + \cdots + d_k) \prod_{1 \leq i \leq k} (d_i+1) .
	\end{eqnarray*}
	}
 We may also compute
\beq\label{gf1}
 \ga(f) := \sum_{0 \leq \al_1 \leq d_1} \cdots \sum_{0 \leq \al_k \leq d_k} (\al_1, \ldots, \al_k) = \frac{1}{2} \Dcal  (d_1,\ldots , d_k).
 \eeq
\xtra{For the density we proceed by induction. As a simple special case in dimension $2$ we see that 
\begin{eqnarray*}
 \sum_{0 \leq \al_1 \leq d_1} \sum_{0 \leq \al_2 \leq d_2}(\al_1, \al_2)
	&=& \sum_{0 \leq \al_1 \leq d_1} ((d_2+1)\al_1, \sum_{0 \leq \al_2 \leq d} \al_2) \\
	&=&  \sum_{0 \leq \al_1 \leq d_1} ((d_2+1)\al_1, \frac{1}{2} d_2(d_2+1)) \\
	&=&  ((d_2+1) \sum_{0 \leq \al_1 \leq d_1} \al_1, \frac{1}{2} d_2(d_2+1)(d_1+1)) \\
	&=&  ((d_2+1)\frac{1}{2}d_1(d_1+1), \frac{1}{2} d_2(d_2+1)(d_1+1)) .
	\end{eqnarray*}
 }

Now we may state our second main result:
\begin{thm}\label{thm_gen_PPW_cor}
For any integer $r \geq R(f)(d+1)$, if $q_i^{\frac{1}{2(r-M(f))}} < H_i < q_i^{\frac{1}{2} + \frac{1}{4(r-M(f))}}$ for each $i=1,\ldots, k$, then
\[S_k(f;\Nbf,\Hbf)
	\ll \| \Hbf\|^{1-\frac{1}{r}} \|\qbf\|^{\frac{-r-M(f)+1}{4r(r-M(f))} + \ep} (\qbf^{\ga(f)})^{\frac{1}{4r(r-M(f))}} q_{\max}^{\frac{2rk}{4r(r-M(f))}}q_{\min}^{-\frac{M(f)}{4r(r-M(f))}} ,
\]
	uniformly in $\Nbf$, with implied constant dependent on $r,d,k,\ep$ and independent of the coefficients of $f$.
	\end{thm}

	\begin{cor}\label{thm_gen_PPW_cor_q}
Under the conditions of Theorem \ref{thm_gen_PPW_cor}, if in addition $q_1 = \cdots = q_k=q$ for a fixed prime $q$, 
\[S_k(f;\Nbf,\Hbf)	\ll  \| \Hbf\|^{1-\frac{1}{r}} q^{\frac{k(r+1-M(f))}{4r(r-M(f))} + \ep}.
\]
	\end{cor}
	
These results rely on a multi-variable version of the Vinogradov Mean Value Theorem tailored to the set $\Fbf(f)$ (see Section \ref{sec_VMT}). Note that the bounds on the right hand side are sharper than those of Theorem \ref{thm_box_Vin} and its corollary  for any $f$ such that the inclusion in (\ref{inclusion}) is strict, so that $M(f) < M_{d,k}$, in which case the bounds also hold for a larger range of $r$ since $R(f) < R_{d,k}$. 

{\bf Example B.} We will highlight the strength of this second type of result by considering the particularly simple case of  
\beq\label{f_example}
f(\xbf) = x_1^{d_1}x_2^{d_2}
\eeq
in dimension $k= 2$ with fixed integers $d_1 > d_2 \geq 1$, and  total degree $d=d_1+d_2$. We compute that
\[ \Fbf (f) = \{ x_1^{\al_1}x_2^{\al_2} : 0 \leq \al_1 \leq d_1, 0 \leq \al_2 \leq d_2 , (\al_1,\al_2) \neq (0,0)\}.\]
Thinking of $d_2$ as fixed and $d_1$ as arbitrarily large relative to $d_2$, we see from (\ref{Rf1})--(\ref{gf1}) that 
\begin{eqnarray*}
 R(f) &=&  (d_1+1)(d_2+1)-1 \approx d_1 \nonumber \\
M(f) &=& \frac{1}{2} (d_1+d_2)(d_1+1)(d_2+1) \approx d_1^2 \label{Mf}\\
 \ga(f) &=& \frac{1}{2}(d_1+1)(d_2+1)(d_1,d_2) \approx (d_1^2,d_1). \nonumber
 \end{eqnarray*}
In comparison, we see from (\ref{RE_dfn_0}) that 
\begin{eqnarray*}
 R_{d,k} &=& \frac{1}{2}(d_1+d_2+2)(d_1+d_2+1) -1 \approx d_1^2 \nonumber \\
M_{d,k} &=& \frac{1}{3}(d_1+d_2+2)(d_1+d_2+1)(d_1+d_2) \approx d_1^3. \label{Md}
\end{eqnarray*}
Thus the bound provided by Corollary \ref{thm_gen_PPW_cor_q}  is significantly sharper than that of Corollary \ref{thm_box_Vin_corq}, and the range $r \geq R(f) (d+1)$ is longer than the range $r \geq R_{d,k}(d+1)$.

This is a genuinely multi-dimensional phenomenon. In the case of dimension $k=1$, given any fixed polynomial $f(x)$ of degree $d$, one necessarily computes 
\[\Fbf(f) =\{ x, x^2, \ldots, x^d\} = \Fbf_{d,1};\]
 that is, in dimension $k=1$, equality always holds in (\ref{inclusion}).  The strength of Theorem \ref{thm_gen_PPW_cor} stems from the fact that in the multi-variable setting, given a fixed polynomial $f$, the resulting set $\Fbf(f)$ is typically much smaller than $\Fbf_{d,k}$.
  
{\bf Example C.} In addition, we note that the explicit presence of the exponent $\ga(f)$ in Theorem \ref{thm_gen_PPW_cor} can also be advantageous, when the primes $q_i$ have varying sizes. (Such a situation can be encountered in applications, for example, which require counting integral points on a hypersurface within a box with disparate side-lengths.) Continuing with the example (\ref{f_example}),  moduli $q_1,q_2$, and degrees $d_1,d_2$ with $d_1$ arbitrarily large relative to $d_2$, the term $\qbf^{\ga(f)}$ in Theorem \ref{thm_gen_PPW_cor} takes the form (for some constants $c_i$)
\[ \qbf^{\ga(f)} \approx q_1^{c_1d_1^2} q_2^{c_2d_1}
;\]
this is advantageous 
compared to the analogous factor in Theorem \ref{thm_box_Vin}, namely
\[ \|\qbf\|^{\frac{M}{k}} \approx q_1^{c_3d_1^3}q_2^{c_4d_1^3}, \]
if for example $q_2$ is large compared to $q_1$.
With these contrasting examples in mind, we now turn to the fully general setting in which we will work for the remainder of the paper.

\section{The general setting of translation-dilation invariant systems}\label{sec_main_thm}

Let $\Fbf$ denote a system of homogeneous polynomials, 
\[ \Fbf = \{F_1,\ldots, F_R\}\]
 with $F_\ell \in \Z[X_1,\ldots, X_k]$ for each $1 \leq \ell \leq R$. Consider for any integer $r \geq 1$ the system of $R$ simultaneous Diophantine equations 
\beq\label{F_sys}
 \sum_{j=1}^r (F_\ell(\xbf_j) - F_\ell(\ybf_j)) = \mathbf{0}, \qquad \text{for all $1 \leq \ell \leq R$,}
 \eeq
where $\xbf_j, \ybf_j \in \Z^k$ for $j = 1,\ldots, r$.
Define $J_r(\Fbf;X)$ to be the number of integral solutions of the system (\ref{F_sys}) with $1 \leq x_{j,i}, y_{j,i} \leq X$ for all $1 \leq j \leq r$, $1 \leq i \leq k$. In \cite{PPW13}, Parsell, Prendiville and Wooley prove strong upper bounds for $J_r(\Fbf;X)$ when $\Fbf$ is a translation-dilation invariant system, which we now define. 

We say $\Fbf$ is a \emph{translation-dilation invariant} system if the following two properties are satisfied: (i) the polynomials $F_1,\ldots, F_R$ are each homogeneous of positive degree; and (ii) there exist polynomials 
\[ c_{m\ell} \in \Z[\xi_1,\ldots, \xi_k], \qquad \text{for each $1 \leq m \leq R, 0 \leq \ell \leq m$},\]
with $c_{mm}=1$ for $1 \leq m \leq R$, such that for any $\boldsymbol{\xi} \in \Z^k$, 
\beq\label{Fc}
 F_m (\xbf + \xibf) = c_{m0}(\xibf) + \sum_{\ell=1}^m c_{m\ell}(\xibf) F_\ell(\xbf), \qquad 1 \leq m \leq R.
 \eeq
(As in \cite{PPW13}, we note that the number of solutions to (\ref{F_sys}) counted by $J_r(\Fbf;X)$ is not affected when one re-orders the $F_\ell$ or takes independent linear combinations of the original forms; so we will say a system is translation-dilation invariant if it is equivalent via such manipulations to a system which is translation-dilation invariant in the strict sense.)

Translation-dilation invariant systems are simple to generate. As a first example, note that $\Fbf_{d,k}$ defined in (\ref{F_dfn}) is a translation-dilation invariant system. As a second example, given any polynomial $f$, the set $\Fbf(f)$ constructed in (\ref{Ff_dfn}) is a translation-dilation invariant system. In fact, more generally,  given any collection of homogeneous polynomials, say 
\[ G_1,\ldots, G_h \in \Z[X_1,\ldots, X_k],\]
one can construct a translation-dilation invariant system.
One first constructs the set $\Gscr$ consisting of all the partial derivatives 
\[ \frac{\partial^{t_1+\cdots + t_k}}{\partial x_1^{t_1} \cdots \partial x_k^{t_k}} G_m(\xbf), \qquad 1 \leq m \leq h,\]
with integral $t_i \geq 0$ for each $1 \leq i \leq k$. The set $\Gscr$ is clearly finite; let $\Gscr_0 = \{ F_1,\ldots, F_R\}$ denote the subset of $\Gscr$ consisting of all polynomials with positive degree, labeled so that $\deg F_1 \leq \deg F_2 \leq \cdots \leq \deg F_R$. Then one confirms via the multi-dimensional Taylor's theorem that the conditions (\ref{Fc}) hold, for some choice of coefficients $c_{m\ell}(\xibf) \in \Z[\xi_1,\ldots, \xi_k]$ such that $c_{mm}(\xibf) = 1$ for $1 \leq m \leq R$.  Furthermore, by replacing the set of forms $\Gscr_0$ by any subset whose span contains $F_1,\ldots, F_R$, we may assume that the set $\{ F_1,\ldots, F_R\}$ is linearly independent, in which case we say the system is \emph{reduced}.
Finally, we introduce the  notion of a \emph{monomial} translation-dilation invariant system, simply by requiring that each form $F_\ell$ in the system be a monomial. We will also avoid certain vacuous cases by making explicit the requirement that a reduced monomial translation-dilation invariant system of dimension $k$ in variables $X_1,\ldots, X_k$ includes for each $i=1,\ldots, k$ at least one monomial of positive degree with respect to $X_i$. To summarize, we may conclude that for any polynomial $f$ we will consider, the set $\Fbf(f)$  is a reduced monomial translation-dilation invariant system.

We now define the parameters used in \cite{PPW13} to characterize a reduced monomial translation-dilation invariant system $\Fbf = \{ F_1,\ldots, F_R\}$ with monomials $F_\ell \in \Z[X_1,\ldots, X_k]$.  We say that $k = k(\Fbf)$ is the \emph{dimension} of the system and $R = R(\Fbf)$ is the \emph{rank}. For each monomial $F_\ell$ we let $d_\ell(\Fbf) = \deg (F_\ell)$ be the total degree of the monomial. We define the \emph{degree} $d = d(\Fbf)$ of the system  by 
\[ d(\Fbf) = \max_{1 \leq \ell \leq R} d_\ell(\Fbf).\]
We define the \emph{weight} $M = M(\Fbf)$ of the system by 
\[ M(\Fbf) = \sum_{\ell=1}^R d_\ell(\Fbf).\]

It is also convenient to use an alternative representation of $\Fbf =\{ F_1,\ldots, F_R \}$ by explicitly writing $\Fbf$ as a collection of monomials 
\[ \{ \xbf^\be : \be \in \Lambda (\Fbf) \},\]
for a fixed collection $\Lambda (\Fbf) $ of $R$ distinct non-zero multi-indices $\be \in \Z_{\geq 0}^k$. If $\Fbf$ has degree $d$, then we see that $|\be| \leq d$ for each $\be \in \Lambda(\Fbf)$ (and there exists some $\be \in \Lambda(\Fbf)$ with $|\be| = d$), and the rank $R(\Fbf)$ is $| \Lambda (\Fbf)|$. The  weight is 
\[ M(\Fbf) = \sum_{\be \in \Lambda (\Fbf)} |\be|.\]
Finally, we define the  notion of the \emph{density} $\ga = \ga(\Fbf) \in \Z_{\geq 0}^k$ of the system by setting 
\beq\label{dens_dfn}
 \ga (\Fbf) = \sum_{\be \in \Lambda (\Fbf)} \be.
 \eeq
In particular, we note that $|\ga| = M(\Fbf)$.

\subsection{Vinogradov Mean Value Theorem}\label{sec_VMT}
We recall the main result of Parsell, Prendiville and Wooley in full generality:
\begin{letterthm}[Theorem 2.1 of \cite{PPW13}]\label{thm_PPW_gen}
Let $\Fbf$ be a reduced translation-dilation invariant system  having dimension $k$, degree $d$, rank $R$ and weight $M$. Suppose that $r$ is a natural number with $r \geq R(d+1)$. Then for each $\ep>0$, 
\beq\label{J_bound_gen}
J_r(\Fbf;X) \ll X^{2rk - M + \ep}.
\eeq
\end{letterthm}
Theorem \ref{thm_PPW} corresponds to the special case of taking $\Fbf$ to be the system $\Fbf_{d,k}$ in (\ref{F_dfn}).

\subsection{Statement of general results}
Our main result in full generality is:
\begin{thm}\label{thm_gen_PPW}
Let $\Fbf$ be a reduced monomial translation-dilation invariant system having dimension $k$, degree $d$, rank $R$, weight $M$, and density $\ga$. Let $\Fscr$ denote the set of all real-valued polynomials spanned by the system $\Fbf$. 
If $r >R(d+1)$ and $ q_i^{\frac{1}{2(r-M)}} < H_i < q_i^{\frac{1}{2} + \frac{1}{4(r-M)}}$ for each $i=1,\ldots, k$, then
\[
\sup_{f \in \Fscr} | S_k(f;\Nbf,\Hbf) |\ll_{r,d,k,\ep}  \| \Hbf\|^{1-\frac{1}{r}}  \|\qbf\|^{\frac{-r-M+1}{4r(r-M)} + \ep} (\qbf^{\ga})^{\frac{1}{4r(r-M)}}q_{\max}^{\frac{2rk}{4r(r-M)}}q_{\min}^{-\frac{M}{4r(r-M)}},
\]
uniformly in $\Nbf $. 
\end{thm}
\begin{cor}\label{thm_gen_PPW_corq}
Under the hypotheses of Theorem \ref{thm_gen_PPW}, if in addition $q_1 = \cdots = q_k = q$ for a fixed prime $q$,
\[
\sup_{f \in \Fscr} | S_k(f;\Nbf,\Hbf) |	\ll_{r,d,k,\ep}  \| \Hbf\|^{1-\frac{1}{r}} q^{\frac{k(r+1-M)}{4r(r-M)} + \ep}.
\]
\end{cor}


As usual, we may check the strength of this result by computing that if $H_i = q^{1/4 + \kappa_i}$ for each $i = 1,\ldots, k$ then Corollary \ref{thm_gen_PPW_corq} provides a bound of size $\| \Hbf \| q^{-\del}$ where 
 \beq\label{del_strength_k}
 \del \approx \frac{ \left( \sum_{i=1}^k \kappa_i \right)^2}{k}.
 \eeq 
 (See Section \ref{sec_thm_Vin_proof} for details.)
Notably, this is independent of the degree, rank, and weight of the system $\Fbf$, and only dependent on the dimension $k$. This also recovers the strength of the original Burgess bound (\ref{Burgess_del}) in dimension $k=1$.  

We note that the input of Theorem \ref{thm_PPW_gen} is crucial; if we used the Burgess method alone without inputting an appropriate Vinogradov Mean Value Theorem, we would obtain a result with 
\beq\label{del_strength_triv}
 \del \approx \frac{\left( \sum_{i=1}^k \kappa_i\right)^2}{M+k}
 \eeq
in place of (\ref{del_strength_k}), which is weaker both because it is smaller and because it is dependent on the degree $d$ of the polynomial $f$. (We will record such a result later in Theorem \ref{thm_gen_triv}.)
Finally, is  clear that Theorems \ref{thm_box_Vin} and  \ref{thm_gen_PPW_cor} are immediate corollaries of Theorem \ref{thm_gen_PPW}. 


%
 
We remark that the approach of this paper is expected to generalize, when suitably adapted, to translation-dilation invariant systems of homogeneous polynomials that are not necessarily monomials. 
Additionally, we note that $\{ x_1,\ldots, x_k\}$ is a special case of a system of $k$ linearly independent linear forms over $\F_q$. In \cite{Bur68} Burgess proved that if $\{L_i\}_{1 \leq i \leq k}$ is a system of $k$ linearly independent linear forms over $\F_q$ for $q$ prime, then 
\beq\label{SHq}
 \sum_{\bstack{\nbf \in \Z^k}{\nbf \in (\Nbf,\Nbf+\Hbf]}} \chi( \prod_{i=1}^k L_i(\nbf)) \ll H^kq^{-\del}
 \eeq
for some small $\del = \del(k) >0$, provided $H>q^{\frac{1}{2} - \frac{1}{2k+2} + \ep}.$
More recently, Chang \cite{Cha09} (for $k=2$) and Bourgain and Chang \cite{BoCh10} (for $k \geq 3$) have proved a bound of the form (\ref{SHq}) that is nontrivial in the original Burgess range of $H>q^{1/4+\ep}$.
It is reasonable to expect that the methods of this paper will generalize to mixed character sums involving products of linear forms of this type.

\subsection{Notation}
For two $k$-tuples $\Kbf = (K_1,\ldots, K_k)$ and $\Hbf = (H_1,\ldots, H_k)$ of real numbers, we will let $\Kbf \leq \Hbf$ represent that all the following conditions hold:
\[ K_1 \leq H_1, \ldots, K_k \leq H_k.\]
We define $\Kbf < \Hbf$ and $\Kbf \ll \Hbf$ similarly.
We will denote by $\Kbf \circ \Hbf$ the coordinate-wise product,
\[ \Kbf \circ \Hbf = (K_1 H_1, \ldots, K_k H_k).\]
We will write 
\[ \Kbf^{-1} = (K_1^{-1}, \ldots, K_k^{-1}),\]
and use the notation
\[  \Hbf/\Kbf = \Hbf \circ \Kbf^{-1} = (H_1/K_1,\ldots, H_k/K_k).\]
For any $k$-tuple $\Kbf = (K_1,\ldots, K_k)$ we set 
\[ \| \Kbf \| = K_1 \cdots K_k.\]
For a scalar $q$, we will say that $\Kbf = (K_1,\ldots, K_k)$ is regarded modulo $q$ if each $K_i$ is regarded modulo $q$. We will say $\Kbf$ is regarded modulo $\Hbf$ if $K_i$ is regarded modulo $H_i$ for each $i=1,\ldots, k$.
For a scalar $q$, we will write $\Kbf q = (K_1 q, \ldots, K_k q)$. We will let implied constants depend on $r,d,k$ and $\ep$ as appropriate. We define the notation $\Lscr(\qbf) = \prod \log q_i$.

\section{Activation of the Burgess method}
Let $\Fbf = \{ F_1, \ldots, F_R\}$ be a given reduced monomial translation-dilation invariant system of dimension $k$, degree $d$, rank $R$, weight $M$ and density $\ga$. 
We will let $\Lambda (\Fbf)$ be the associated set of multi-indices, so that we can represent $\Fbf$ as $\{\xbf^\be : \be \in \Lambda(\Fbf)\}$. We will let $\Fscr(\Fbf)$ denote the set of all real-valued polynomials spanned by the set of monomials comprising $\Fbf$. We will let $\Fscr_0(\Fbf)$ denote the set of all real-valued polynomials spanned by $1 \union \Fbf$; that is, we expand $\Fscr(\Fbf)$ to include polynomials with constant terms. We correspondingly set $\Lambda_0(\Fbf) = \{ (0,\ldots, 0)\} \union \Lambda (\Fbf)$.

The family $\Fscr_0(\Fbf)$ is invariant under translations: by the relations (\ref{Fc}), if $f(\xbf) \in \Fscr_0(\Fbf)$ then $f(\xbf + \xibf) \in \Fscr_0(\Fbf)$ for all $\xibf \in \R^k$. 
Similarly $\Fscr_0(\Fbf)$ is invariant under dilations $\xbf \mapsto \xibf \circ \xbf$: that is, if $f(\xbf) \in \Fscr_0(\Fbf)$ then $f(\xibf \circ \xbf) \in \Fscr_0(\Fbf)$ for all $\xibf \in \R^k$. This is a stronger type of dilation invariance than dilation by scalars, and is a consequence of using monomial systems. To confirm this, we simply represent $f$ as 
\[ f(\xbf) = \sum_{\be \in \Lambda_0 (\Fbf)} f_\be \xbf^\be \]
with coefficients $f_\be$, so that 
\[ f(\xibf \circ \xbf) = \sum_{\be \in \Lambda_0 (\Fbf)} f_\be (\xibf \circ \xbf)^\be =  \sum_{\be \in \Lambda_0 (\Fbf)} (f_\be \xibf^\be)  \xbf^\be, \]
which is also a polynomial in $\Fscr_0(\Fbf)$. Finally, we note that since we assume in the definition of a reduced monomial translation-dilation invariant system that for each $i=1,\ldots, k$, $\Fbf$ contains a monomial of positive degree in $X_i$, expanding the relations (\ref{Fc}) using the multinomial theorem shows that  linear monomials in each of $X_1,\ldots, X_k$ also belong to $\Fbf$.
We will use these facts repeatedly in the argument to come.

From now on $\Fbf$ will be the fixed system given above. Fix primes $q_1,\ldots, q_k$ (not necessarily distinct) and let $\qbf=(q_1,\ldots, q_k)$. For each $i=1,\ldots,k$ let $\chi_i$ be a non-principal multiplicative Dirichlet character modulo $q_i$.
Instead of working directly with $S_k(f;\Nbf,\Hbf)$ we will define 
\[ T(\Fbf; \Nbf,\Hbf) = \sup_{f \in \Fscr_0(\Fbf)} \sup_{\Kbf \leq \Hbf} \left| \sum_{\bstack{\xbf \in \Z^k}{\xbf \in (\Nbf, \Nbf+\Kbf]}} e(f(\xbf)) \chi_1(x_1)\cdots \chi_k(x_k) \right|,\]
which certainly majorizes $S_k(f;\Nbf,\Hbf)$.
We first note that $T(\Fbf; \Nbf, \Hbf)$ is unchanged if the supremum over $f \in \Fscr_0 (\Fbf)$ is restricted to $f \in \Fscr (\Fbf)$, as appears in the statement of our theorems. Second, we note that $T(\Fbf;\Nbf;\Hbf)$ is periodic modulo $\qbf$ with respect to $\Nbf$. Indeed, if $\Nbf = \Mbf \circ \qbf + \Lbf$ for an integer tuple $\Mbf$, we can express $ T(\Fbf; \Nbf,\Hbf)$ as
 \begin{multline*}
 \sup_{f \in \Fscr_0(\Fbf)} \sup_{\Kbf \leq \Hbf} \left| \sum_{\bstack{\xbf \in \Z^k}{\xbf \in (\Lbf, \Lbf+\Kbf]}} e(f(\xbf + \Mbf \circ \qbf)) \chi_1(x_1 + M_1 q_1)\cdots \chi_k(x_k+M_kq_k) \right| \\
 = \sup_{f \in \Fscr_0(\Fbf)} \sup_{\Kbf \leq \Hbf} \left| \sum_{\bstack{\xbf \in \Z^k}{\xbf \in (\Lbf, \Lbf+\Kbf]}} e(f(\xbf + \Mbf \circ \qbf)) \chi_1(x_1 )\cdots \chi_k(x_k) \right| ,
 \end{multline*}
which is $T(\Fbf; \Lbf,\Hbf),$ as claimed.
Thus we see that it suffices to consider $\Nbf$ with $0 \leq N_i < q_i$ for $i=1,\ldots, k$.
We also note that in $T(\Fbf;\Nbf,\Hbf)$ it suffices to regard the coefficients of the polynomial $f$ modulo $1$; by a compactness argument, one sees that the value of $T(\Fbf;\Nbf,\Hbf)$ is achieved by a particular choice of polynomial $f$ and length $\Kbf$.
 
We now begin the familiar opening gambit of the Burgess method. Given a fixed $\Hbf = (H_1,\ldots, H_k)$, we let $P_1, \ldots P_k$ be a set of parameters each satisfying $1 \leq P_i \leq H_i$, to be chosen precisely later. For each $i=1,\ldots, k$ we fix a set of primes
 \[ \Pcal_i= \{ P_i < p \leq 2P_i \}.\]
We then let $\Pscr$ denote the corresponding set of $k$-tuples of primes:
 \[ \Pscr = \{ \pbf = (p_1,\ldots, p_k) : p_i \in \Pcal_i, i=1,\ldots, k \}.\]
Since we will restrict to $H_i =o(q_i)$ in our theorems, we will be able to assume $p_i \ndiv q_i$ for all $p_i \in \Pcal_i$, for all $i$. We also note that for each $i$, $|\Pcal_i| \gg P_i (\log P_i)^{-1} \gg P_i (\log q_i)^{-1}$, so that 
\beq\label{Psize}
 |\Pscr| \gg P_1 \cdots P_k (\prod_{i=1}^k \log q_i)^{-1}  = \| \Pbf \| \Lscr (\qbf)^{-1}.
\eeq

Fix a tuple $\Kbf \leq \Hbf$ and a tuple $\pbf$ of primes in $\Pscr$; then each $\xbf \in (\Nbf, \Nbf + \Kbf]$ may be split into residue classes modulo $\pbf$, so that for each $i=1,\ldots, k$, we may write
\[ x_i = a_i q_i + p_i m_i,\]
where $0 \leq a_i < p_i$ and $m_i  \in (N_i^{a_i,p_i},N_i^{a_i,p_i}+K_i^{a_i,p_i}],$ where we have set
\begin{eqnarray*}
N_i^{a_i,p_i} &=& \frac{N_i-a_i q_i}{p_i}, \\
K_i^{a_i,p_i}  &= &\frac{K_i}{p_i} \leq \frac{H_i}{p_i} \leq \frac{H_i}{P_i}.
\end{eqnarray*}
 That is to say, $\xbf = \abf \circ \qbf + \pbf \circ \mbf$ with $\boldsymbol{0} \leq \abf < \pbf$ and $\mbf \in [\Nbf^{\abf, \pbf}, \Nbf^{\abf, \pbf} + \Kbf^{\abf,\pbf})$.
Then we see that 
\begin{multline}\label{efchi}
 \sum_{\bstack{\xbf \in \Z^k}{\xbf \in (\Nbf, \Nbf+\Kbf]}} e(f(\xbf)) \chi_1(x_1)\cdots \chi_k(x_k)\\
	 = \sum_{\boldsymbol{0} \leq \abf < \pbf} \sum_{\mbf \in (\Nbf^{\abf,\pbf}, \Nbf^{\abf,\pbf} + \Kbf^{\abf,\pbf}]} e(f(\abf \circ \qbf +\pbf \circ \mbf)) \prod_{i=1}^k \chi_i (a_iq_i + p_im_i).
	 \end{multline}
We  may remove the dependence of the multiplicative characters on $p_i$, since 
\[ \prod_{i=1}^k \chi_i (a_iq_i + p_im_i) = \prod_{i=1}^k \chi_i (m_i) \cdot \prod_{i=1}^k \chi_i (p_i).\]
	Thus after taking absolute values and taking the supremum over $f \in \Fscr_0(\Fbf)$ and $\Kbf \leq \Hbf$ in (\ref{efchi}), we see that 
	\[ T(\Fbf;\Nbf,\Hbf) \leq \sum_{\boldsymbol{0} \leq \abf < \pbf} T(\Fbf; \Nbf^{\abf,\pbf}, \Hbf /\Pbf).\]
After averaging over the set $\Pscr$, we then have 
\beq\label{TTap}
 T(\Fbf;\Nbf,\Hbf) \leq |\Pscr|^{-1} \sum_{\pbf \in \Pscr} \sum_{\boldsymbol{0} \leq \abf < \pbf} T(\Fbf; \Nbf^{\abf,\pbf}, \Hbf / \Pbf).
 \eeq

We will now make the starting points $\Nbf^{\abf,\pbf}$ of the sums $T(\Fbf; \Nbf^{\abf,\pbf}, \Hbf / \Pbf)$ independent of $\abf, \pbf$ via the following lemma:
\begin{lemma}\label{lemma_TUL}
For any tuple $\Ubf$ of real numbers and tuple $\Lbf$ of  real numbers with $L_i \geq 1$ for $i=1,\ldots, k$,
\[ T(\Fbf; \Ubf,\Lbf)  \leq  2^{2k}  \| \Lbf \|^{-1} \sum_{\Ubf - \Lbf < \mbf \leq \Ubf} T(\Fbf; \mbf,2\Lbf).\]
\end{lemma}
Suppose $T(\Fbf; \Ubf, \Lbf)$ is attained by a certain polynomial $f$ and a tuple $\Kbf \leq \Lbf$; then we write
\[ T(\Fbf; \Ubf,\Lbf) = \left| \sum_{\xbf \in (\Ubf, \Ubf + \Kbf]} e(f(\xbf)) \prod_{i=1}^k \chi_i(x_i) \right| .\]
By the inclusion-exclusion principle, for any fixed $\Rbf$ with $\Ubf - \Lbf < \Rbf \leq \Ubf$,
 \[ \sum_{\xbf \in (\Ubf, \Ubf + \Kbf]} e(f(\xbf)) \prod_{i=1}^k \chi_i(x_i)  =  \sum_{\bstack{{\boldsymbol{\ep}} = (\ep_1, \ldots, \ep_k)}{\ep_i \in \{ 0,1\}}} (-1)^{|\boldsymbol{\ep}|} \sum_{\xbf \in (\Rbf,( \mathbf{1} - \boldsymbol{\ep}) \circ \Kbf + \Ubf]} e(f(\xbf)) \prod_{i=1}^k \chi_i(x_i) . \]
Here $|\boldsymbol{\ep}| = \ep_1 + \cdots + \ep_k$, and $\mathbf{1} - \boldsymbol{\ep} = (1-\ep_1, \ldots, 1 - \ep_k)$.

We next note that for any $\Rbf$ with $\Ubf - \Lbf < \Rbf \leq \Ubf$ and any $\boldsymbol{\ep}$ as above, the side-lengths of the box $(\Rbf, (\mathbf{1}-\boldsymbol{\ep}) \circ \Kbf + \Ubf]$ satisfy
\[ (\mathbf{1} - \boldsymbol{\ep}) \circ \Kbf+ \Ubf - \Rbf \leq 2 \Lbf.
\]
Thus 
\[ \left | \sum_{\xbf \in (\Ubf, \Ubf + \Kbf]} e(f(\xbf)) \prod_{i=1}^k \chi_i(x_i)   \right|
	\leq 2^k T(\Fbf; \Rbf, 2\Lbf).\]
	We finally note that there are at least $L_i/2$ integers in the interval $(U_i - L_i, U_i]$ and hence at least $2^{-k}\|\Lbf\|$ choices for tuples $\Rbf$ in the box $(\Ubf - \Lbf, \Ubf]$, so that averaging over all of these choices produces the result of Lemma \ref{lemma_TUL}.
	
We now apply Lemma \ref{lemma_TUL} to (\ref{TTap}) with the choice $\Lbf = \Hbf/\Pbf$, to see by (\ref{Psize}) that 
\begin{eqnarray*}
 T(\Fbf;\Nbf,\Hbf)& \leq& 2^{2k}\| \Hbf/\Pbf \|^{-1} 
 	 |\Pscr|^{-1}  \sum_{\pbf \in \Pscr} \sum_{\boldsymbol{0} \leq \abf < \pbf} 
	 	\sum_{\Nbf^{\abf,\pbf} - \Hbf/\Pbf < \mbf \leq \Nbf^{\abf,\pbf}} T(\Fbf; \mbf, 2\Hbf/\Pbf) \\
		& \ll & \| \Hbf \|^{-1}  \Lscr(\qbf) \sum_{\pbf \in \Pscr} \sum_{\boldsymbol{0} \leq \abf < \pbf} 
	 	\sum_{\Nbf^{\abf,\pbf} - \Hbf/\Pbf < \mbf \leq \Nbf^{\abf,\pbf}} T(\Fbf; \mbf, 2\Hbf/\Pbf) .
		\end{eqnarray*}
Now for each $\mbf$ we define $\Acal(\mbf)$ to be the quantity 
 \[  \# \{ \abf, \pbf : 0 \leq a_i < p_i \; \text{and} \; p_i \in \Pcal_i: \frac{N_i - a_iq_i}{p_i} - \frac{H_i}{P_i} < m_i 	\leq  \frac{N_i - a_iq_i}{p_i}  , i=1,\ldots, k \}.\]
With this notation, we may now write
\[  T(\Fbf;\Nbf,\Hbf) \ll \| \Hbf \|^{-1} \Lscr (\qbf)  \sum_{\mbf \in \Z^k} \Acal (\mbf) T(\Fbf; \mbf, 2\Hbf/\Pbf).\]
We now define 
\[ S_1 = \sum_{\mbf} \Acal(\mbf),\]
and 
\[ S_2 = \sum_{\mbf} \Acal(\mbf)^2.\]
We record the following facts, which we prove in Section \ref{sec_A_lemma}:
\begin{lemma}\label{lemma_A}
We have $\Acal(\mbf)=0$ unless $|m_i| \leq 2q_i$ for $i=1,\ldots, k$. Furthermore if $H_iP_i<q_i$ for each $i=1,\ldots, k$ then 
\[ S_1 \leq S_2 \ll \| \Hbf \| \, \| \Pbf \|.\]
\end{lemma}

After a repeated application of H\"{o}lder's inequality, Lemma \ref{lemma_A} allows us to conclude that 
\begin{eqnarray*}
T(\Fbf;\Nbf,\Hbf) &\ll & \| \Hbf \|^{-1} \Lscr (\qbf) S_1^{1 - \frac{1}{r}} S_2^{\frac{1}{2r}} \left\{ \sum_{\bstack{\mbf}{|m_i| \leq 2q_i}} T(\Fbf;\mbf, 2\Hbf/\Pbf)^{2r} \right\}^{\frac{1}{2r}} \\
& \ll & \| \Hbf\|^{-\frac{1}{2r}} \| \Pbf \|^{1 - \frac{1}{2r}} \Lscr (\qbf) \left\{ \sum_{\bstack{\mbf}{|m_i| \leq 2q_i}} T(\Fbf; \mbf, 2\Hbf/\Pbf)^{2r} \right\}^{\frac{1}{2r}}.
	\end{eqnarray*}
We now recall that $T(\Fbf; \mbf,\Kbf)$ is periodic in $\mbf$ with respect to $\qbf$, so that it suffices to write
\beq\label{TT2r}
 T(\Fbf;\Nbf,\Hbf) \ll  \| \Hbf\|^{-\frac{1}{2r}} \| \Pbf \|^{1 - \frac{1}{2r}} \Lscr (\qbf) \left\{ \sum_{\mbf \modd{ \qbf}} T(\Fbf; \mbf, 2\Hbf/\Pbf)^{2r} \right\}^{\frac{1}{2r}}.
 \eeq

We now make the step of removing the supremum over lengths in the definition of $T(\Fbf;\Nbf,\Kbf)$. We define for any tuples $\Mbf,\Kbf$ with $K_i >0$ the sum
\[ T_0 (\Fbf;\Mbf,\Kbf) = \sup_{f \in \Fscr_0(\Fbf)} \left| \sum_{\Mbf < \xbf \leq \Mbf + \Kbf} e(f(\xbf)) \chi_1(x_1) \cdots \chi_k(x_k) \right|.\]
We will use the following lemma, a $k$-dimensional version of Lemma 2.2 of Bombieri and Iwaniec \cite{BoIw86}, whose proof we indicate in Section \ref{sec_BoIw}.
\begin{lemma}\label{lemma_BoIw}
Let $a(\nbf)$ be a sequence of complex numbers indexed by tuples $\nbf$ supported on the set $\nbf \in (\Abf,\Abf+\Bbf]\subset \Z^k.$ Let $I = (\Cbf,\Cbf + \Dbf]$ be any product of intervals with $I \subseteq (\Abf, \Abf+\Bbf]$. Then 
\[ \sum_{\nbf \in I} a(\nbf) \ll ( \prod_{i=1}^k \log (B_i +2)) \sup_{\theta \in \R^k} \left| \sum_{\nbf \in (\Abf, \Abf + \Bbf]} a(\nbf) e(\theta \cdot \nbf) \right|.\]
\end{lemma}
This lemma allows us to relate $T(\Fbf; \Mbf,\Kbf)$ to $T_0(\Fbf; \Mbf,\Kbf)$ since as long as $d \geq 1$, Lemma \ref{lemma_BoIw} shows that
\[ T(\Fbf; \Mbf,\Kbf) \ll ( \prod_{i=1}^k \log (K_i +2)) T_0(\Fbf; \Mbf,\Kbf) .\]
Note that here we use the assumption that $d\geq1$, so that the linear exponential factor accrued in the application of Lemma \ref{lemma_BoIw} is absorbed in the supremum over polynomials $f \in \Fscr_0(\Fbf)$.
We also henceforward assume that $K_i <q_i$ for each $i=1,\ldots, k$, so that the logarithmic factor is bounded above by $\ll \Lscr(\qbf)$; this condition will be satisfied by our final choice of $K_i$, as we will later verify.

We may now re-write (\ref{TT2r}) as 
\beq\label{TS3}
 T(\Fbf;\Nbf,\Hbf) \ll  \| \Hbf\|^{-\frac{1}{2r}} \| \Pbf \|^{1 - \frac{1}{2r}} \Lscr(\qbf)^2 S_3(2\Hbf/\Pbf) ^{\frac{1}{2r}},
 \eeq
where we define
\beq\label{S3_dfn}
S_3(\Kbf) :=  \sum_{\mbf \modd{\qbf}} T_0(\Fbf; \mbf, \Kbf )^{2r}.
 \eeq

\section{Approximation of polynomials}
We will now bound $S_3(\Kbf)$, focusing first on an individual sum $T_0(\Fbf; \mbf,\Kbf)$; recall that we assume from now on that each $K_i < q_i$. 
As in \cite{HBP14a}, the key step is to remove the supremum over all polynomials in $\Fscr_0(\Fbf)$ by showing, roughly speaking, that two polynomials with coefficients that are sufficiently close may be regarded as producing equivalent contributions, and thus we will majorize the supremum by summing over a collection of representative polynomials. We first perform a dissection of the coefficient space of $\Fscr_0(\Fbf)$, recalling that we may regard the coefficients of any $f \in \Fscr_0(\Fbf)$ modulo $1$. 

We recall the collection of multi-indices $\Lambda_0 (\Fbf) = \{ (0,\ldots, 0)\} \union \Lambda(\Fbf)$ associated to the system $\Fbf$. Since $\Fbf$ has rank $R$, we have $R=|\Lambda (\Fbf)|$ and $R+1=|\Lambda_0(\Fbf)|$, so that $R+1$ is the dimension of the coefficient space of $\Fscr_0(\Fbf)$. 

  Fix positive integers $Q_1,\ldots, Q_k$ and set $\Qbf = (Q_1,\ldots, Q_k)$. We will choose $Q_i$ precisely later; for now we assume that $Q_i \geq K_i$ for each $i$, which we will verify later. 
  We index  the coefficient space $[0,1]^{R+1}$ as 
 \[ [0,1]^{R+1} = [0,1] \times \cdots \times[0,1]= \prod_{\be \in \Lambda_0 (\Fbf)} [0,1]^{(\be)} .\]
 For each of the $R+1$ multi-indices $\be \in \Lambda_0(\Fbf)$, we partition the corresponding unit interval $[0,1]^{(\be)} = [0,1]$ indexed by $\be$ into $\Qbf^\be = Q_1^{\be_1}\cdots Q_k^{\be_k}$ sub-intervals of length $(\Qbf^\be)^{-1}$. We claim this partitions the full space $[0,1]^{R+1}$ into $\Qbf^\ga$ boxes, where we recall that $\ga = \ga(\Fbf)$ is the density of the system $\Fbf$, as defined in (\ref{dens_dfn}). 
 We may verify this as follows:  clearly the number of boxes is
 \[  \prod_{\be \in \Lambda_0 (\Fbf)} \Qbf^{\be} = \Qbf^{\del},\]
say, where we have defined
\[ \del = \sum_{\be \in \Lambda_0 (\Fbf)} \be = \sum_{\be \in \Lambda(\Fbf)} \be.\]
This last expression is precisely the definition of the density $\ga = \ga (\Fbf)$. 

We will denote this dissection of the coefficient space as a union 
\beq\label{box_diss}
 [0,1]^{R+1}= \bigcup_\al B_\al 
 \eeq
over $\Qbf^\ga$ many boxes $B_\al$; we may think of $\al$ as a parameter in $\Z_{\geq 0}$ indexing over a fixed ordering of the boxes. We will also associate to each box $B_\al$ the fixed tuple $\theta_\al \in B_\al$ that is the vertex of $B_\al$ with the least value in each coordinate. Thus if we have fixed some enumeration $\be^{(0)}, \ldots, \be^{(R)}$ of the $R+1$ distinct multi-indices $\be \in \Lambda_0 (\Fbf)$, the distinguished vertex $\theta_\al$ of a box $B_\al$ takes the form 
\beq\label{theta_spec}
 \theta_\al = (\theta_{\al,\be^{(0)}}, \ldots, \theta_{\al, \be^{(R)}}) = (c_{\be^{(0)}}\Qbf^{-\be^{(0)}} , \ldots, c_{\be^{(R)}}\Qbf^{-\be^{(R)}}),
 \eeq
where for each $j=0, \ldots, R$, $c_{\be_j}$ is an integer with $0 \leq c_{\be_j} \leq \Qbf^{\be_j} - 1$.
Finally, for any fixed $\theta \in [0,1]^{R+1}$, we define an associated real-valued polynomial on $\R^k$ by
\beq\label{theta_poly_dfn}
 \theta(\Xbf) := \sum_{\be \in \Lambda_0 (\Fbf) } \theta_\be \Xbf^\be.
 \eeq
We note that for any $\theta \in [0,1]^{R+1}$ this polynomial belongs to $\Fscr_0(\Fbf)$.

For any tuple $\mbf$ of integers and any tuple $\taubf$ of positive real numbers and a fixed index $\al$ of a box $B_\al$ with associated vertex $\theta_\al$, we define
 \beq\label{TalF}
  T(\al, \Fbf;\mbf,\taubf) := \left| \sum_{\mathbf{0} < \nbf \leq \taubf} e(\theta_\al(\nbf)) \chi_1(n_1+m_1) \cdots \chi_k(n_k+m_k) \right|.
  \eeq
 Roughly speaking, our goal is to show that for any $\mbf,\Kbf$ there exists a suitable $\al$ such that $T_0(\Fbf; \mbf,\Kbf)$ (which takes a supremum over $f \in \Fscr_0(\Fbf)$) is well approximated by $T(\al, \Fbf; \mbf,\Kbf)$ (which corresponds to the single polynomial $\theta_\al(\Xbf)$). In order to do so, we must use summation by parts, for which we require some notation.

Given any partition $I \union J$ of the set of indices $\{1,\ldots, k\}$ and a $k$-tuple $\nbf$, we will let $\nbf_{(I)}$ denote the tuple of $n_j$ with $j \in I$ and similarly $\nbf_{(J)}$ the tuple of $n_j$ with $j \in J$; thus for example we may write $\nbf = (\nbf_{(I)}, \nbf_{(J)})$ (with some abuse of notation with respect to ordering).
Given a sequence $a(\nbf)$ of complex numbers indexed by $\nbf \in \N^k$, we will define partial summation of $a(\nbf)$ with respect to such partitions as follows:
\[A_{(I),(J)}(\tbf_{(I)}, \sbf_{(J)}) := \sum_{\bstack{0<n_j \leq t_j}{j \in I}} \sum_{\bstack{0<n_j \leq s_j}{j \in J}} a(\nbf).\]

  More specifically, in our application, given a partition $I \union J$ of $\{1,\ldots, k\}$,
a tuple $\mbf$ of integers and tuples $\sbf,\tbf \in \R^k$ of positive real numbers and a fixed index $\al$ of a box  $B_\al$, we define
 \[ T_{(I),(J)}(\al, \Fbf;\mbf, \tbf_{(I)},\sbf_{(J)}) := \left| \sum_{\bstack{0< n_j \leq t_j}{j \in I}} \sum_{\bstack{0< n_j \leq s_j}{j \in J}} e(\theta_\al(\nbf)) \left( \prod_{i=1}^{k} \chi_i(n_i+m_i) \right) \right|.\]

 The key approximation lemma is as follows:
 \begin{lemma}\label{lemma_T_part_sum}
 Given integral tuples $\mbf$ and $\Qbf$ with
 \beq\label{Q_big}
 Q_i \geq K_i \quad \text{for each $i=1,\ldots, k$},
 \eeq
the above dissection provides an index $\al$ such that 
 \beq\label{T0T}
T_0(\Fbf; \mbf,\Kbf) \ll_{k,d} \sum_{J \subseteq \{1,\ldots, k\}} \left( \prod_{j \in J} K_j^{-1} \right)
 \idotsint_{\bstack{(0,K_j]}{ j \in J}}  T_{(^cJ),(J)}(\al, \Fbf;\mbf, \Kbf_{(^cJ)},\tbf_{(J)}) d\tbf_{(J)}.
  \eeq
 Here the sum is over all subsets $J \subseteq \{1,\ldots, k\}$, with corresponding complement $^c J = \{ 1,\ldots, k\} \setminus J$. If $J = \{ j_1,\ldots, j_v\}$ then we set $ d \tbf_{(J)}= dt_{j_1} \cdots dt_{j_v}$.
 \end{lemma}
To prove this, we first observe that for an integral tuple $\mbf$, 
 \begin{eqnarray*}
  T_0 (\Fbf; \mbf,\Kbf) &=& \sup_{f \in \Fscr_0(\Fbf)} \left| \sum_{\mbf < \xbf \leq \mbf + \Kbf} e(f(\xbf)) \chi_1(x_1) \cdots \chi_k(x_k) \right| \\
  	& = &  \sup_{f \in \Fscr_0 (\Fbf)} \left| \sum_{\mathbf{0}< \xbf \leq \Kbf} e(f(\xbf)) \chi_1(x_1+m_1) \cdots \chi_k(x_k+m_k) \right| .
	\end{eqnarray*}
We now write
\[ T_0(\Fbf; \mbf, \Kbf) =  \left| \sum_{\mathbf{0}< \xbf \leq \Kbf} e(f(\xbf)) \chi_1(x_1+m_1) \cdots \chi_k(x_k+m_k) \right|,\]
for some fixed polynomial $f \in \Fscr_0(\Fbf)$, which we write explicitly as 
\[ f(\Xbf) = \sum_{\be \in \Lambda_0(\Fbf)} f_\be \Xbf^\be,\]
where as before we may assume that each $f_\be \in [0,1]$. Given our dissection of the coefficient space $[0,1]^{R+1}$, we may choose a box $B_\al$ with  index $\al$ and distinguished vertex $\theta_\al = (\theta_{\al,\be})_{\be}$ such that  
\beq\label{del_small}
 |f_\be - \theta_{\al,\be}| \leq \Qbf^{-\be}, \quad \text{for each multi-index $\be \in \Lambda_0 (\Fbf)$}.
 \eeq
 (This is simply choosing $\al$ such that the coefficient tuple $(f_\be)_\be$ lies in the box $B_\al$.)
For notational convenience, we will temporarily set $\del_\be = f_\be - \theta_{\al,\be}$ for each  $\be \in \Lambda_0 (\Fbf)$. We then write
\begin{multline}\label{sum_ab_apply}
  \sum_{\mathbf{0}< \xbf \leq \Kbf} e(f(\xbf)) \chi_1(x_1+m_1) \cdots \chi_k(x_k+m_k)\\
	 =  \sum_{\mathbf{0}< \xbf \leq \Kbf} e\left( \sum_{\be \in \Lambda_0 (\Fbf)} \del_\be \xbf^\be \right) e(\theta_\al(\xbf)) \chi_1(x_1+m_1) \cdots \chi_k(x_k+m_k).
	 \end{multline}

 We now apply summation by parts, in the following form, which we prove in Section \ref{sec_part_sum_proof}: 
\begin{lemma}\label{lemma_part_sum_gen}
Let $a(\nbf)$ be a sequence of complex numbers indexed by $\nbf \in \Z^k \intersect (\boldsymbol{0},\Nbf]$.
Let $b(\xbf)$ be a $C^{(k)}$ function on $\R^k$ such that there is a tuple $\Bbf = (B_1,\ldots, B_k)$ of positive real numbers such that for every multi-index 
\beq\label{ga_form}
\kappa = (\kappa_1,\ldots, \kappa_k) \quad \text{with $\kappa_i \in \{0,1\}$} 
\eeq
 we have 
\beq\label{b_assp}
\left|  \frac{\partial^{|\kappa|}}{\partial x_1^{\kappa_1} \cdots \partial x_k^{\kappa_k}} b(\xbf) \right| \leq B_1^{\kappa_1} \cdots B_k^{\kappa_k} = \Bbf^{\kappa} \quad \text{for all $\xbf \in (\mathbf{0},\Nbf]$.}
 \eeq
Then 
\[ \left| \sum_{\nbf\leq \Nbf} a(\nbf) b(\nbf) \right| \leq \sum_{J \subseteq \{1,\ldots, k\}} \left( \prod_{j \in J} B_j \right)
 \idotsint_{\bstack{(0,N_j]}{ j \in J}} \left| A_{(^cJ),(J)}(\Nbf_{(^cJ)},\tbf_{(J)}) \right| d\tbf_{(J)}.\]
 Here the sum is over all subsets $J \subseteq \{1,\ldots, k\}$, with corresponding complement $^c J = \{ 1,\ldots, k\} \setminus J$. If $J = \{ j_1,\ldots, j_v\}$ then we set $ d \tbf_{(J)}= dt_{j_1} \cdots dt_{j_v}$.\end{lemma}
(Note that if $\ka=(0,0,\ldots,0)$ then (\ref{b_assp}) is simply the assumption that $|b(\xbf)| \leq 1$.)
We apply this lemma to (\ref{sum_ab_apply}) with the choices $\Nbf = \Kbf$ and 
\begin{eqnarray*}
 a(\xbf) &=& e(\theta_\al(\xbf)) \chi_1(x_1+m_1) \cdots \chi_k(x_k+m_k),\\
 b(\xbf) & = & e\left( \sum_{\be \in \Lambda_0(\Fbf)} \del_\be \xbf^\be \right) .
 \end{eqnarray*}
We may verify that for a fixed index $j$, if we let $e_j = (0, \ldots, 1, \ldots,0)$ be the $j$-th unit multi-index, then for $\xbf \in (\boldsymbol{0},\Kbf]$,
\begin{eqnarray*}
 \left| \frac{\partial}{\partial x_j} e \left(  \sum_{\be \in \Lambda_0(\Fbf)} \del_\be \xbf^\be \right)  \right|
& = &  \left| \frac{\partial}{\partial x_j} \left( 2\pi i \sum_{\be \in \Lambda_0(\Fbf)} \del_\be \xbf^\be \right)  \right| \\
& \leq & 2\pi  \sum_{\bstack{\be \in \Lambda_0(\Fbf)}{\be_j \geq 1}} \be_j |\del_\be| |\xbf^{\be-e_j}| \\
& \leq & 2\pi  \sum_{\bstack{\be \in \Lambda_0(\Fbf)}{\be_j \geq 1}} \be_j |\del_\be| |\Kbf^{\be-e_j}| \\
& \leq & 2\pi  \sum_{\bstack{\be \in \Lambda_0(\Fbf)}{\be_j \geq 1}} \be_j \Qbf^{-\be} \Kbf^{\be-e_j} \\
& \ll_{k,d}& K_j^{-1},
\end{eqnarray*}
where we have used the assumption (\ref{del_small}) on the size of $\del_\be$, followed by the assumption (\ref{Q_big}) that $Q_i \geq K_i$.

Similarly, one may compute that for each fixed $\kappa$ of the form (\ref{ga_form}), 
\begin{eqnarray*}
 \left| \frac{\partial^{|\kappa|}}{\partial x_1^{\kappa_1} \cdots \partial x_k^{\kappa_k}} e \left(  \sum_{\be \in \Lambda_0(\Fbf)} \del_\be \xbf^\be \right)  \right|
& \ll_{k,d} \Kbf^{-\kappa},
\end{eqnarray*}
so that a bound of the form (\ref{b_assp}) is satisfied with $B_i = K_i^{-1}$.
We may thus apply Lemma \ref{lemma_part_sum_gen} to (\ref{sum_ab_apply}) to conclude that (\ref{T0T}) holds.

A repeated application of H\"{o}lder's inequality in (\ref{T0T}) shows that $ T_0(\Fbf; \mbf, \Kbf)^{2r}$ is at most
\[   \ll_{k,r,d} \sum_{J \subseteq \{1,\ldots, k\}} \left( \prod_{j \in J} K_j^{-1} \right) \idotsint_{\bstack{(0,K_j]}{ j \in J}}T_{(^cJ),(J)}(\al, \Fbf;\mbf, \Kbf_{(^cJ)}, \tbf_{(J)})^{2r} d\tbf_{(J)}.\]
\xtra{We apply H\"{o}lder's to write
\begin{eqnarray*}
|T_0|^{2r} & \ll& (2^k)^{2r-1} \sum_{J} \left( \idotsint (\prod K_j^{-1}) T dt \right)^{2r} \\
	& \ll_{k,r}& \sum_{J} \left( \idotsint_{(0,K_j]} (\prod K_j^{-1})^{\frac{2r}{2r-1}} \right)^{2r-1} \idotsint |T|^{2r} \\
	& \ll& \sum_J (\prod K_j^{-1}) \idotsint |T|^{2r}.
\end{eqnarray*}
}
This is still for the fixed index $\al$ provided by Lemma \ref{lemma_T_part_sum}; as in \cite{HBP14a}, in order to eliminate the awkward dependence on $\al$, we sum trivially on the right hand side over all values of the parameter $\al$ that indexes the boxes in the dissection (\ref{box_diss}), so that  by positivity, $T_0(\Fbf; \mbf,\Kbf)^{2r}$ is at most
\begin{multline*}
    \ll_{k,r,d} \sum_{J \subseteq \{1,\ldots, k\}} \left( \prod_{j \in J} K_j^{-1} \right) \idotsint_{\bstack{(0,K_j]}{ j \in J}} \sum_{\al} T_{(^cJ),(J)}(\al, \Fbf;\mbf, \Kbf_{(^cJ)}, \tbf_{(J)})^{2r} d\tbf_{(J)} .
   \end{multline*}
Now we sum over $\mbf \modd{\qbf}$, so that $ \sum_{\mbf \modd{\qbf}} T_0(\Fbf; \mbf,\Kbf)^{2r}$ is at most
\[ 
\ll  \sum_{J \subseteq \{1,\ldots, k\}} \left( \prod_{j \in J} K_j^{-1} \right) \idotsint_{\bstack{(0,K_j]}{ j \in J}} \sum_{\al} \sum_{\mbf \modd{\qbf}} T_{(^cJ),(J)}(\al, \Fbf;\mbf, \Kbf_{(^cJ)}, \tbf_{(J)})^{2r} d\tbf_{(J)} .
\]

We note that if $\tbf \leq \Kbf$ then  by positivity, for any index set $J$,
\beq\label{TTsup}
   \sum_{\al} \sum_{\mbf \modd{\qbf}}  T_{(^cJ),(J)}(\al, \Fbf;\mbf, \Kbf_{(^cJ)},\tbf_{(J)})^{2r} \leq \sup_{\taubf \leq \Kbf} \left( \sum_{\al}   \sum_{\mbf \modd{\qbf}}T(\al,\Fbf;\mbf, \taubf)^{2r} \right),
\eeq
where $T(\al,\Fbf;\mbf,\taubf)$ is defined by (\ref{TalF}).  Applying this in the integrand above and noting the normalization complementing the region of integration, we see that 
\[  \sum_{\mbf \modd{\qbf}}T_0(\Fbf; \mbf,\Kbf)^{2r} \ll_{k,r,d} \sup_{\taubf \leq \Kbf} \left( \sum_{\al}  \sum_{\mbf \modd{\qbf}}T(\al,\Fbf;\mbf, \taubf)^{2r} \right).\]
For convenience, we define
\[ S_4(\taubf) = \sum_\al \sum_{\mbf \modd{\qbf}} T(\al, \Fbf;\mbf, \taubf)^{2r}.\]
We may conclude:
\begin{lemma}\label{lemma_TS4}
\beq\label{T0S4}
S_3(\Kbf)  = \sum_{\mbf \modd{\qbf}} T_0(\Fbf; \mbf,\Kbf)^{2r} \ll_{k,r,d} \sup_{\taubf \leq \Kbf} S_4(\taubf).
\eeq
\end{lemma}

We now summarize what we have proved so far, by combining the result of Lemma \ref{lemma_TS4} with (\ref{TS3}) and (\ref{S3_dfn}):
\begin{prop}\label{prop_TS4}
As long as $\Qbf \geq \Kbf$, $\Kbf = 2\Hbf/\Pbf < \qbf$, and $H_iP_i< q_i$ for all $i=1,\ldots, k$,
\[
 T(\Fbf;\Nbf,\Hbf) \ll  \| \Hbf\|^{-\frac{1}{2r}} \| \Pbf \|^{1 - \frac{1}{2r}} \Lscr (\qbf)^2 \left \{ \sup_{\taubf \leq 2\Hbf/\Pbf} S_4 (\taubf) \right\}^{\frac{1}{2r}}.
\]
\end{prop}
 Thus our focus  turns to bounding $S_4(\taubf)$ for any fixed tuple $\taubf$ with $\taubf \leq \Kbf= 2\Hbf/\Pbf$. 
We recall the definition of the boxes $B_\al$, the vertex $\theta_\al$ associated to each box $B_\al$, and the associated polynomial $\theta_\al(\xbf)$ defined by (\ref{theta_poly_dfn}). We will represent a set of cardinality $2r$ of tuples $\xbf^{(j)}  = (x^{(j)}_1, \ldots, x^{(j)}_k ) \in \Z^k$ by $\{ \xbf \} = \{\xbf^{(1)}, \ldots, \xbf^{(2r)}\}$.
We set 
\[ \Sig_A (\{ \xbf\}) := \sum_\al e \left( \sum_{j=1}^{2r} \ep(j) \theta_\al (\xbf^{(j)}) \right),\]
where $\ep(j) := (-1)^j$. We note that we may trivially bound $\Sig_A$ by the number of summands, that is the number of boxes, namely 
\beq\label{Sig_A_triv}
|\Sig_A (\{ \xbf\}) | \leq \Qbf^\ga,
\eeq
where $\ga = \ga(\Fbf)$ is the density of the system $\Fbf$.

For each $i=1,\ldots, k$ let $\Del_i(q_i)$ denote the order of $\chi_i$ modulo $q_i$; furthermore for each $j=1,\ldots, 2r$ set
$\del_i(j)=1$ if $j$ is even and $\del_i(j)=\Delta_i(q_i)-1$ if $j$ is odd.
Now define for each $i=1,\ldots,k$  the single-variable polynomial 
\[ G_i(\Del_i(q_i), \{ \xbf \};X) := \prod_{j=1}^{2r} (X+ x^{(j)}_i)^{\del_i(j)}.\]
Finally, set
\[ \Sig_B(\{ \xbf \}; \qbf) :=  \prod_{i=1}^k \left( \sum_{m_i = 1}^{q_i} \chi_i( G_i(\Del_i(q_i), \{ \xbf \};m_i))\right).\]
 
We now expand the sums in the definition of $S_4 (\taubf)$ to see that with this notation,
\beq\label{S4sigsig}
 S_4(\taubf) =  \sum_{\bstack{\xbf^{(1)},\ldots, \xbf^{(2r)} \in \Z^k}{ \mathbf{0} < \xbf^{(j)} \leq \taubf}} \Sig_A(\{\xbf\}) \Sig_B (\{ \xbf\};\qbf).
	\eeq
We will now proceed in two parts: first, we will sum trivially over the boxes $B_\al$, using the trivial bound (\ref{Sig_A_triv}) for $\Sig_A$. This will result in the following proposition:
\begin{prop}\label{prop_S_box_triv}
Under the assumption that the tuple $\Kbf$ satisfies 
\[ q_i^{\frac{1}{2r}} \ll K_i \ll q_i^{\frac{1}{2r}} \qquad \text{for each $i=1,\ldots, k$},
\]
 we have
\[\sup_{\taubf \leq \Kbf} S_4(\taubf) \ll  \Qbf^\ga \| \Kbf \|^{2r}  \| \qbf\|^{\frac{1}{2}},\]
 where $\ga$ is the density of the system $\Fbf$ and the implied constant depends on $r,d,k$.
\end{prop}

As an immediate consequence we will prove:
\begin{thm}\label{thm_gen_triv}
Let $\Fbf$ be a reduced monomial translation-dilation invariant system  having dimension $k$, degree $d$, rank $R$, weight $M$ and density $\ga$. Let $\Fscr$ denote the set of all real-valued polynomials spanned by the system $\Fbf$. 
If $r \geq 1$ and $q_i^{\frac{1}{2r}} < H_i < q_i^{\frac{1}{2} + \frac{1}{4r}}$ for each $i=1,\ldots, k$, then
\[
\sup_{f \in \Fscr} | S_k(f;\Nbf,\Hbf) |\ll_{r,d,k}  \| \Hbf\|^{1-\frac{1}{r}}  \|\qbf \|^{\frac{r+1}{4r^2}} (\qbf^{\ga})^{\frac{1}{4r^2}}\Lscr(\qbf)^{2}.
\]
\end{thm}

This is the result that leads to (\ref{del_strength_triv}) in the case $q_1 = \cdots = q_k =q$.  
We will improve on this  in Section \ref{sec_add} by proving a nontrivial upper bound for $\Sig_A$ via Theorem \ref{thm_PPW_gen},  which we will  consequently use to give a refinement (Proposition \ref{prop_S_box_nontriv}) of Proposition \ref{prop_S_box_triv}.  Our main general result, Theorem \ref{thm_gen_PPW}, will  then follow from this refinement. For purposes of comparison, we note that when for example $q_1 = \cdots = q_k=q$, Theorem \ref{thm_gen_PPW} is sharper than Theorem \ref{thm_gen_triv} for $r > M+k$, as well as in the sense that  (\ref{del_strength_k}) improves on (\ref{del_strength_triv}).

 \section{The multiplicative component}
 Our treatment of the multiplicative component $\Sig_B$ is the same for both Theorem \ref{thm_gen_triv} and Theorem \ref{thm_gen_PPW}, and hinges upon an application of the Weil bound.
 It will be convenient to regard a collection $\{ \xbf \}$ as either a set of cardinality $2r$ of $k$-tuples $\xbf^{(1)}, \ldots, \xbf^{(2r)} \in \Z^k$, or equivalently as a set of cardinality $k$ of $2r$-tuples, which we will denote by $\zbf^{(1)}, \ldots, \zbf^{(k)} \in \Z^{2r}$; in matrix form we are regarding the $\xbf^{(j)}$ as the rows of a $2r \times k$ matrix, and the $\zbf^{(i)}$ as the columns. We will denote such a collection $\zbf^{(1)},\ldots, \zbf^{(k)}$ also by $\{ \zbf \}$.
 
 We recall the definitions of $\Del_i(q_i)$ and $\del_i(j)$. 
 We now define for any $2r$-tuple $\zbf = (z_1,\ldots, z_{2r})$ the single-variable polynomial
 \[ F(\Del_i(q_i), \zbf; X):= \prod_{j=1}^{2r} (X+ z_j)^{\del_i(j)}.\]
 We can now write
\beq\label{Sig_z}
 \Sig_B(\{ \xbf \}; \qbf) = \Sig_B(\{ \zbf \};\qbf) = \prod_{i=1}^k \Sig_B(\zbf^{(i)}; \chi_i,q_i) 
 \eeq
where we define
\[  \Sig_B(\zbf^{(i)}; \chi_i,q_i) :=  \sum_{m_i = 1}^{q_i} \chi_i\left( F(\Del_i(q_i),\zbf^{(i)};m_i)\right).\]
We  aim to apply the following consequence of the well-known Weil bound:
\begin{lemma}\label{lemma_Weil}
Let $\chi$ be a character of order $\Delta(q) >1$ modulo a prime $q$. Suppose that $F(X)$ is a polynomial which is not a perfect $\Del(q)$-th power in $\overline{\mathbb{F}}_q[X]$. Then 
\[ \left| \sum_{m=1}^q \chi(F(m)) \right| \leq (\deg (F)-1)\sqrt{q}.\]
\end{lemma}
For a fixed $i$, we can apply Lemma \ref{lemma_Weil} to bound 
\beq\label{Sig_B_good}
 \Sig_B(\zbf^{(i)}; \chi_i,q_i) \ll q_i^{1/2},
 \eeq
  unless $F(\Del_i(q_i),\zbf^{(i)};X)$ is a perfect $\Del_i(q_i)$-th power over $\overline{\mathbb{F}}_{q_i}$, in which case
we apply the trivial bound
  \beq\label{Sig_B_bad}
 \Sig_B(\zbf^{(i)}; \chi_i,q_i) \ll q_i.
 \eeq
 
   We define a $2r$-tuple $\zbf = (z_1,\ldots, z_{2r})$ to be \emph{bad} if for each $j=1,\ldots, 2r$ there exists $\ell \neq j$ such that $z_\ell =z_j$, and to be \emph{good} otherwise. 
We have the following simple statement:
\begin{lemma}\label{lemma_poly_bad}
Fix a character $\chi$ with order $\Del(q)>1$ modulo a prime $q$. 
Fix a tuple $\zbf = (z_1,\ldots, z_{2r})$ with $0 < z_j \leq u$ for each $j=1,\ldots, 2r$. If $u \leq q$ and $F(\Del(q),\zbf;X)$ is a perfect $\Del(q)$-th power modulo $q$, then $\zbf$ is bad. 
\end{lemma}
This result is clear, since if $\zbf$ were good, there would be a value $y$, say, which is taken only by $z_j$ for one index $j \in \{ 1, \ldots, 2r\}$; thus the factor $(X+y)$ would appear with multiplicity 1 or $\Del(q)-1$ in $F(\Del(q), \zbf;X)$, neither of which is divisible by $\Del(q)$.

Lemma \ref{lemma_poly_bad} is useful for a single factor $\Sig_B(\zbf^{(i)};\chi_i,q_i)$, but we must also consider how many tuples in a collection $\{\zbf\}$ are bad. 
For each subset $\Scal \subseteq \{1,\ldots, k\}$ (possibly empty), we say a collection $\{\zbf\} = \{\zbf^{(1)},\ldots, \zbf^{(k)}\}$ is $\Scal$-bad if $\zbf^{(i)}$ is bad for $i \in \Scal$ and good for $i \not\in \Scal$. 
For each subset $\Scal \subseteq \{ 1,\ldots, k\}$,  we let $\Bcal(\Scal;\taubf)$ denote the set of collections $\{\zbf\}$ that are $\Scal$-bad and such that  for each $1 \leq i \leq k$ the entries of $\zbf^{(i)}$ are at most $\tau_i$. 
(Implicitly this also specifies when the original tuple $\{\xbf\}$ belongs to $\Bcal(\Scal;\taubf)$; we will use this in the computation (\ref{S4N}).) 
We now prove an upper bound on the cardinality of the set $\Bcal(\Scal;\taubf)$:
\begin{lemma}\label{lemma_bad_count} 
For any fixed subset $\Scal \subseteq \{ 1,\ldots, k\}$,
\[ \#\Bcal (\Scal;\taubf) \leq r^{k(2r+1)}   \left( \prod_{i \in \Scal} \tau_i^{r} \right)\left( \prod_{i \not\in \Scal} \tau_i^{2r} \right) .
 \]
\end{lemma}
 We recall from the classical Burgess method (see for example Lemma 3.2 of \cite{HBP14a}) that there are at most $r^{2r+1}u^r$ choices for a single bad $2r$-tuple $\zbf$ with entries at most $u$.
 Fix a subset $\Scal \subseteq\{ 1,\ldots, k\}$. For each $i \in \Scal$ there are then at most $r^{2r+1}\tau_i^r$ choices for $\zbf^{(i)}$ that are bad, and for each $i \not \in \Scal$ there are trivially at most $\tau_i^{2r}$ ways to choose $\zbf^{(i)}$ that are good. 
Thus upon recalling $|\Scal| \leq k$, the lemma is proved.

 Finally, we note:
 \begin{lemma}\label{lemma_sig_bad}
Suppose $\{\zbf \} \in \Bcal(\Scal,\taubf)$ and let $\sig(\Scal)$ be the indicator multi-index for $\Scal$, that is $\sig (\Scal) = (\sig_1,\ldots, \sig_k)$  with $\sig_i = 1$ if $i \in \Scal$ and $\sig_i = 0$ if $i \not\in\Scal$.  If $\taubf$ is such that $\tau_i \leq q_i$ for each $i=1,\ldots, k$, then
 \beq\label{sig_bounded}
 \Sig_B(\{\zbf\};\qbf) \ll_{r,k} \|\qbf\|^{\frac{1}{2}} \qbf^{\sig(\Scal)/2}.
 \eeq
 \end{lemma} 
 We simply note that within the product (\ref{Sig_z}) we may apply the Weil bound (\ref{Sig_B_good})  for each $ i \not\in \Scal$ and the trivial bound (\ref{Sig_B_bad})  for each $i \in \Scal$; this suffices for the lemma.

 We now consider (\ref{S4sigsig}), applying the trivial bound (\ref{Sig_A_triv}) to $\Sig_A$ and decomposing $\Sig_B$ as follows:
 \begin{eqnarray*}
 S_4(\taubf) 
	&=& \sum_{\bstack{\xbf^{(1)},\ldots ,\xbf^{(2r)} \in \Z^k}{ \mathbf{0} < \xbf^{(j)} \leq \taubf}} \Sig_A(\{\xbf\}) \Sig_B (\{ \xbf\};\qbf) \\
	& \leq & \Qbf^\ga	 \sum_{\bstack{\xbf^{(1)},\ldots ,\xbf^{(2r)} \in \Z^k}{ \mathbf{0} < \xbf^{(j)} \leq \taubf}} \left| \Sig_B (\{ \xbf\};\qbf) 	\right|			\\
	& = & \Qbf^\ga	\sum_{\Scal \subseteq \{1, \ldots, k\}  }	\sum_{ \{ \zbf \} \in \Bcal(\Scal;\taubf)} \left|	\Sig_B (\{ \zbf\};\qbf) 	\right|			\\
	& \ll_{r,k}&\Qbf^\ga \sum_{\Scal \subseteq \{1, \ldots, k\}  } \left( \prod_{i \in \Scal} \tau_i^{r} \right)\left( \prod_{i \not\in \Scal} \tau_i^{2r} \right)   \| \qbf\|^{1/2} \qbf^{\sig(\Scal)/2}.
 \end{eqnarray*}
Here we have applied Lemmas \ref{lemma_bad_count} and \ref{lemma_sig_bad}.
We now note that since this is monotone in each $\tau_i$, we have the upper bound
\[ \sup_{\taubf \leq \Kbf} S_4(\taubf) \ll \Qbf^\ga \sum_{\Scal \subseteq \{1, \ldots, k\}  }\left( \prod_{i \in \Scal} K_i^{r} \right)\left( \prod_{i \not\in \Scal} K_i^{2r} \right)   \| \qbf\|^{1/2} \qbf^{\sig(\Scal)/2}.
\]
We re-write this as 
\[ \sup_{\taubf \leq \Kbf} S_4(\taubf)  \ll \Qbf^\ga\|\Kbf\|^{2r} \| \qbf \|^{1/2}\left\{ 1 +  \sum_{\bstack{ \Scal \subseteq \{1,\ldots, k\}}{\Scal \neq \emptyset}} \Kbf^{-r \sig(\Scal)}\qbf^{\sig(\Scal)/2}\right\}.
\]

Under the assumption
 \beq\label{tau_assp}
q_i^{\frac{1}{2r}} \ll  K_i \ll q_i^{\frac{1}{2r}}
\quad \text{for each $i=1,\ldots, k$},
\eeq
we have $\Kbf^{-r\sig(\Scal)} \qbf^{\sig(\Scal)/2} = O(1)$ for each subset $\Scal$ and as a result
 \[ \sup_{\taubf \leq \Kbf} S_4(\taubf) \ll_{r,k}  \Qbf^\ga \| \Kbf \|^{2r}  \| \qbf \|^{\frac{1}{2}} ,\]
 which proves Proposition \ref{prop_S_box_triv}.
 
 \subsection{Proof of Theorem \ref{thm_gen_triv}}
With Proposition \ref{prop_S_box_triv} in hand, it is simple to derive Theorem \ref{thm_gen_triv}.
Recalling that $\Kbf = 2\Hbf/\Pbf$, the condition (\ref{tau_assp}) leads us to 
choose the parameters $P_i$ such that
\[ \frac{1}{2} H_i q_i^{- \frac{1}{2r}} \leq P_i < H_i q_i^{-\frac{1}{2r}}, \quad \text{for each $i=1,\ldots, k$} ,\]
in which case (\ref{tau_assp}) holds. We also note that as long as $H_i > q_i^{\frac{1}{2r}}$, we may choose $P_i \geq 1$. We furthermore choose $Q_i = \lceil K_i \rceil$ for each $i$. 
 In order to satisfy the further conditions $H_i P_i < q_i$ of Lemma \ref{lemma_A}, we must restrict to ranges with $H_i < q_i^{\frac{1}{2} + \frac{1}{4r}}$.

 With these choices, 
 we input the result of Proposition \ref{prop_S_box_triv} into Proposition \ref{prop_TS4} to see that 
\begin{eqnarray*}
 T(\Fbf;\Nbf,\Hbf) &\ll&  \| \Hbf\|^{-\frac{1}{2r}} \| \Pbf \|^{1 - \frac{1}{2r}} \Lscr(\qbf)^2 \left \{  \Qbf^\ga \| \Hbf \|^{2r} \| \Pbf \|^{-2r}  \| \qbf \|^{\frac{1}{2}}  \right\}^{\frac{1}{2r}} \\
  & \ll&  \| \Hbf\|^{1-\frac{1}{2r}} \| \Pbf \|^{- \frac{1}{2r}}  (\Hbf/\Pbf)^{\frac{\ga}{2r}} \|\qbf \|^{\frac{1}{4r}}  \Lscr(\qbf)^2 .
  \end{eqnarray*}
 We now note that because of our choice of $\Pbf$, 
 \[ (\Hbf/\Pbf)^\frac{\ga}{2r} \ll (\qbf^\ga)^{\frac{1}{4r^2}}.\]
Thus we may conclude that
 \begin{eqnarray*} 
  T(\Fbf;\Nbf,\Hbf) & \ll&   \| \Hbf\|^{1-\frac{1}{r}} \|\qbf\|^{  \frac{1}{4r^2}} (\qbf^\ga)^{\frac{1}{4r^2}} \| \qbf\|^{\frac{1}{4r}}  \Lscr(\qbf)^2  \\
& \ll &  \| \Hbf\|^{1-\frac{1}{r}} \| \qbf\|^{\frac{r+1}{4r^2}} (\qbf^\ga)^{\frac{1}{4r^2}} \Lscr(\qbf)^2  ,
  \end{eqnarray*}
which proves Theorem \ref{thm_gen_triv}.
In particular, if $q_j=q$ for all $j$, this simplifies to
\[  T(\Fbf;\Nbf,\Hbf)
 \ll    \| \Hbf\|^{1-\frac{1}{r}} q^{\frac{k(r+1) +M}{4r^2}} (\log q)^{2k},  \]
 where we recall that $M = M(\Fbf)$ is the weight of the system $\Fbf$.

\subsection{Optimal choice of $r$}\label{sec_thm_triv_proof}
We make a remark on the case $q_1=\cdots = q_k =q$ and the optimal choice of $r$. Suppose that for each $i$, $H_i = q^{1/4 + \kappa_i}$. Set $\sig = \sum_{i=1}^k \kappa_i.$ Then Theorem \ref{thm_gen_triv} provides an upper bound of the size $\| H \| q^{-\del + \ep}$ where 
$ \del = (\sig r- \frac{1}{4}(k + M))r^{-2}$.
As a function of $r$, this attains a maximum at the real value $ r_0 = (k + M)(2 \sig )^{-1}.$
Choosing $r = r_0 + \theta$ where $-1/2 \leq \theta < 1/2$, we see that as claimed in (\ref{del_strength_triv}), $\del$ is approximately of size $ \del \approx \sig^2(M+k)^{-1}.$

\section{The additive component: nontrivial analysis}\label{sec_add}
 We now return to a nontrivial analysis of the additive component $\Sig_A$, which will lead to our main result Theorem \ref{thm_gen_PPW}. Our goal is to connect the analysis of $\Sig_A$ to  a Vinogradov Mean Value Theorem for the translation-dilation invariant system $\Fbf$.
 We again recall the definition of the boxes $B_\al$ that partition the coefficient space $[0,1]^{R+1}$, and in particular the definition (\ref{theta_spec}) of the distinguished vertex $\theta_\al$ associated to each box $B_\al$. It is convenient to recall the fixed ordering $\be^{(0)}, \ldots, \be^{(R)}$ of the multi-indices $\be \in \Lambda_0(\Fbf)$. We compute:
 \begin{eqnarray*}
 \Sig_A (\{ \xbf\}) &=& \sum_\al e \left( \sum_{j=1}^{2r} \ep(j) \theta_\al (\xbf^{(j)}) \right) \\
 & = & \sum_\al e \left( \sum_{\be \in \Lambda_0(\Fbf)} \theta_{\al,\be} \left( \sum_{j=1}^{2r} \ep(j)( \xbf^{(j)})^\be \right) \right) \\
  & = & \sum_\al e \left( \sum_{\be = \be^{(0)}, \ldots, \be^{(R)} } \theta_{\al,\be} \left( \sum_{j=1}^{2r} \ep(j)( \xbf^{(j)})^\be \right) \right) \\
  & = &\sum_{c_{\be^{(0)}}, \ldots, c_{\be^{(R)}}}  e \left( \sum_{\be = \be^{(0)}, \ldots, \be^{(R)}} c_\be \Qbf^{-\be} \left( \sum_{j=1}^{2r} \ep(j)( \xbf^{(j)})^\be \right) \right) ,
  \end{eqnarray*}
 where the sum over $c_{\be^{(0)}}, \ldots, c_{\be^{(R)}}$ indicates summing  for each $i=0,\ldots, R$ the parameter $c_{\be^{(i)}}$ over integers $0 \leq c_{\be^{(i)}} \leq \Qbf^{\be^{(i)}}-1$.
 Thus
 \begin{eqnarray*}
 \Sig_A (\{ \xbf\})
  & = &\sum_{c_{\be^{(0)}}, \ldots, c_{\be^{(R)}}}  \left\{ \prod_{\be = \be^{(0)}, \ldots, \be^{(R)}} e \left(c_\be \Qbf^{-\be} \left( \sum_{j=1}^{2r} \ep(j)( \xbf^{(j)})^\be \right) \right) \right\} \\
 & = & \prod_{\be = \be^{(0)}, \ldots, \be^{(R)}}  \left\{\sum_{c_\be \modd{\Qbf^\be}} e \left(c_\be \Qbf^{-\be} \left( \sum_{j=1}^{2r} \ep(j)( \xbf^{(j)})^\be \right) \right) \right\}\\
  & = & \prod_{\be \in \Lambda_0(\Fbf)}  \left\{\sum_{c_\be \modd{\Qbf^\be}} e \left(c_\be \Qbf^{-\be} \left( \sum_{j=1}^{2r} \ep(j)( \xbf^{(j)})^\be \right) \right) \right\}.
 \end{eqnarray*}
Since the multi-index $\be=(0,\ldots,0)$ contributes only a multiplicative factor of 1, we see that
\[  \Sig_A (\{ \xbf\})= \prod_{\be \in \Lambda(\Fbf)}  \left\{\sum_{c_\be \modd{\Qbf^\be}} e \left(c_\be \Qbf^{-\be} \left( \sum_{j=1}^{2r} \ep(j)( \xbf^{(j)})^\be \right) \right) \right\}.\]
By orthogonality of characters we therefore have
\[  \Sig_A (\{ \xbf\})= \Qbf^\ga \;  \Xi_\Qbf (\Fbf;\{ \xbf \})\]
 where 
 $\Xi_\Qbf(\Fbf; \{\xbf \})$ is the indicator function for the set
 \beq\label{xi_set_Q}
  \{ \xbf^{(1)},\ldots,  \xbf^{(2r)} \in \mathbb{Z}^{k} \intersect (\mathbf{0},\taubf ] : \sum_{j=1}^{2r} \ep(j)( \xbf^{(j)})^\be \con 0 \modd{\Qbf^\be}, \forall  \be \in \Lambda(\Fbf) \}.
  \eeq
 Here we have also used the fact, previously observed, that 
 \[  \prod_{\be \in \Lambda(\Fbf)} \Qbf^\be =  \Qbf^\ga.\]
 In our application we will have $\tau_i \leq K_i= 2H_i/P_i$ for each $i=1,\ldots, k$. So far we have only assumed that $Q_i \geq K_i$ for each $i$; we now furthermore assume that each $Q_i$ is sufficiently large that the congruences in the definition of the set (\ref{xi_set_Q}) must be identities in $\Z$. We check that for any multi-index $\be \in \Lambda(\Fbf)$ and any collection $\{ \xbf \}$ in the set (\ref{xi_set_Q}), 
 \[ \left| \sum_{j=1}^{2r} \ep(j)( \xbf^{(j)})^\be \right| < 2r \taubf^\be\leq 2r \Kbf^\be \leq (2r \Kbf)^\be  = (4r\Hbf/\Pbf)^\be .\]
 Thus we choose 
 \beq\label{Q_choice}
  Q_i  = \lceil 4r H_i/P_i \rceil \quad \text{for each $i=1,\ldots, k$}.
  \eeq
 With this choice the congruences in (\ref{xi_set_Q}) must be identities in $\Z$, and we may replace $\Xi_\Qbf(\Fbf; \{\xbf\})$ by the indicator function $\Xi (\Fbf; \{ \xbf \})$ of the set 
 \[ V_{r} (\Fbf; \taubf) :=  \{ \xbf^{(1)},\ldots,  \xbf^{(2r)} \in \mathbb{Z}^{k} \intersect (\mathbf{0},\taubf ] : \sum_{j=1}^{2r} \ep(j)( \xbf^{(j)})^\be =0, \forall \be \in \Lambda(\Fbf)\}.\]
  We have shown:
 \begin{prop}\label{prop_VJT}
 Given a collection $\{ \xbf \} = \{ \xbf^{(1)}, \ldots, \xbf^{(2r)}\}$ with $\xbf^{(j)} \in \Z^k \intersect (\mathbf{0},\taubf]$ for each $j = 1,\ldots, 2r$,
if $\taubf \leq \Kbf  = 2\Hbf/\Pbf$ and we choose $\Qbf$ as in (\ref{Q_choice}), then
  \[ \Sig_A (\{ \xbf\})  = \Qbf^\ga \, \Xi (\Fbf; \{ \xbf \}).\]
 \end{prop}
 
 We now set 
 \[ \tau_{\max} = \max \{ \tau_1,\ldots, \tau_k\},\]
and note that  
 \[ \#V_r(\Fbf;\taubf) \leq J_{r}(\Fbf; \tau_{\max}),\]
 where $J_{r}(\Fbf; X)$ is the counting function for the system of equations (\ref{F_sys}) corresponding to the given reduced monomial translation-dilation invariant system $\Fbf$. We recall from (\ref{J_bound_gen}) that the conjectured upper bound is 
 \beq\label{J_bound_gen'}
 J_{r} (\Fbf;X) \ll X^{2rk-M +\ep},
 \eeq
 known unconditionally for $r \geq R(d+1)$ by Theorem \ref{thm_PPW_gen} (Section \ref{sec_VMT}).

 We now return to the consideration of $S_4(\taubf)$, given in (\ref{S4sigsig}) in terms of the additive component $\Sig_A$ and the multiplicative component $\Sig_B$. 
Define for each subset $\Scal \subseteq \{1,\ldots, k \}$ the quantity
  \[ N(\Scal; \taubf) = \# \{ \Bcal(\Scal;\taubf) \intersect V_{r} (\Fbf; \taubf) \}.\]
 With $\taubf \leq \Kbf$ and $\Qbf$ as above, we apply Proposition \ref{prop_VJT} and Lemma \ref{lemma_sig_bad} to see that
 \begin{eqnarray}
  S_4(\taubf)& = &\sum_{\bstack{\xbf^{(1)},\ldots, \xbf^{(2r)} \in \Z^k}{ \mathbf{0} < \xbf^{(j)} \leq \taubf}} \Sig_A(\{\xbf\}) \Sig_B (\{ \xbf\};\qbf) \nonumber \\
  & = & \Qbf^\ga \sum_{\bstack{\xbf^{(1)},\ldots, \xbf^{(2r)}\in \Z^k}{ \mathbf{0} < \xbf^{(j)} \leq \taubf}} \Xi (\Fbf; \{ \xbf \})\Sig_B (\{ \xbf\};\qbf) \nonumber \\
  & = &\Qbf^\ga   \sum_{\Scal \subseteq \{ 1,\ldots,k \}} \sum_{\{ \xbf \} \in \Bcal(\Scal; \taubf)} \Xi (\Fbf; \{ \xbf \})\Sig_B (\{ \xbf\};\qbf) \nonumber \\
  & \ll & \Qbf^\ga \|\qbf\|^{1/2}  \sum_{\Scal \subseteq \{ 1,\ldots,k \}} N(\Scal; \taubf)  \qbf^{\sig(\Scal)/2}.
  	\label{S4N}
  \end{eqnarray}

We now bound $N(\Scal; \taubf)$ for each subset $\Scal$.
If $\Scal = \emptyset$, the size of $\Bcal(\Scal;\taubf)$ is no smaller than $O(\tau_{\max}^{2kr})$, so the key restriction comes from the fact that we are counting points that also lie in $V_{r}(\Fbf; \taubf)$. Thus for $\Scal = \emptyset$ we use the estimate
 \[ N(\Scal;\taubf) \leq \# V_{r} (\Fbf; \taubf) \leq J_{r}(\Fbf; \taumax) \ll \taumax^{2rk -M + \ep}, \]
under the assumption that the bound (\ref{J_bound_gen'}) for $J_{r}(\Fbf; X)$ holds. 
For any non-empty $\Scal$, we use an upper bound based only on the size of $\Bcal(\Scal; \taubf)$ from Lemma \ref{lemma_bad_count}: 
 \[ N(\Scal;\taubf) \leq \# \Bcal(\Scal;\taubf) \ll \left( \prod_{i \in S} \tau_i^r \right) \left( \prod_{i \not\in S} \tau_i^{2r} \right) .\]
 We then have 
\[   S_4(\taubf) \ll 
    \Qbf^\ga \left\{ \taumax^{2rk - M+\ep} \|\qbf\|^{\frac{1}{2}} +  \|\qbf \|^{\frac{1}{2}}\sum_{\bstack{\Scal \subseteq \{ 1,\ldots,k \}}{\Scal \neq \emptyset}} \left( \prod_{i\in \Scal} \tau_i^{r} \right) \left( \prod_{i \not\in \Scal} \tau_i^{2r} \right) \qbf^{\sig(\Scal)/2}\right\} .\]
    
    We henceforward assume that 
    \beq
    r>M,
    \eeq 
    so that certainly $2rk-M \geq 0$.
We now define $K_{\max} = \max \{ K_1,\ldots, K_k\}$ and $K_{\min} = \min \{ K_1,\ldots, K_k\}$ and use the fact that the above bound for $S_4(\taubf)$ is monotone in $\taubf \leq \Kbf $. Then
$\sup_{\taubf \leq \Kbf}   S_4(\taubf)$ is at most
\begin{eqnarray*}
& \ll &    \Qbf^\ga \left\{ K_{\max}^{2rk - M+\ep} \|\qbf\|^{\frac{1}{2}} +  \|\qbf \|^{\frac{1}{2}}\sum_{\bstack{\Scal \subseteq \{ 1,\ldots,k \}}{\Scal \neq \emptyset}} \left( \prod_{i\in \Scal}K_i^{r} \right) \left( \prod_{i \not\in \Scal} K_i^{2r} \right) \qbf^{\sig(\Scal)/2}\right\} \\
	& \ll &   \Qbf^\ga    \left\{ K_{\max}^{2rk - M+\ep} \| \qbf \|^{\frac{1}{2}} + \|\Kbf \|^{2r} \| \qbf \|^{\frac{1}{2}}\sum_{\bstack{\Scal \subseteq \{ 1,\ldots,k \}}{\Scal \neq \emptyset}} \Kbf^{-r \sig(\Scal)} \qbf^{\sig(\Scal)/2} \right\}\\
     &  = & \Qbf^\ga    \left\{ K_{\max}^{2rk - M+\ep} \| \qbf \|^{\frac{1}{2}} + \|\Kbf \|^{2r} \| \qbf \|^{\frac{1}{2}}\sum_{\bstack{\Scal \subseteq \{ 1,\ldots,k \}}{\Scal \neq \emptyset}} \left(  \prod_{i \in \Scal} K_i^{-r}q_i^{1/2} \right) \right\} .
     \end{eqnarray*}
     
  Now we make the assumption that for every $i =1,\ldots, k$ we have $K_i^{-r}q_i^{1/2} \leq 1$, or equivalently
  \beq\label{Ki_qi}
   K_i \geq q_i^{\frac{1}{2r}},
  \eeq
 so that the largest contribution from the sum over subsets $\Scal$ comes from sets of cardinality one. 
 Then we have 
 \[
 \sup_{\taubf \leq \Kbf}   S_4(\taubf)
 	\ll
	 \Qbf^\ga   \left\{ K_{\max}^{2rk - M+\ep} \| \qbf \|^{\frac{1}{2}} + \|\Kbf \|^{2r} \| \qbf \|^{\frac{1}{2}} \sum_{i=1}^kK_i^{-r}q_i^{1/2} \right\}.
	\]
	This implies
 \[
 \sup_{\taubf \leq \Kbf}   S_4(\taubf)
 	\ll
	 \Qbf^\ga   \left\{ K_{\max}^{2rk+\ep}K_{\min}^{-M} \| \qbf \|^{\frac{1}{2}} + K_{\max}^{2rk} \| \qbf \|^{\frac{1}{2}} \sum_{i=1}^kK_i^{-r}q_i^{1/2} \right\}.
	\]	
The first term in braces dominates all other terms as long as for each $i=1,\ldots, k$ we have 
\beq\label{Kmin}
 K_{\min}^{-M} \gg K_i^{-r}q_i^{1/2},
 \eeq
which is certainly implied by the condition
\beq\label{Kmin'}
K_i^{r -M} \geq q_i^{1/2};
\eeq
we note that this condition would also guarantee (\ref{Ki_qi}).

 We have proved:
 \begin{prop}\label{prop_S_box_nontriv}
 If $r> M$ and (\ref{Kmin'}) holds for each $i=1,\ldots, k$, then
\[
 \sup_{\taubf \leq \Kbf} S_4(\taubf) \ll \Qbf^\ga K_{\max}^{2rk+\ep}K_{\min}^{ - M} \| \qbf \|^{\frac{1}{2}}.
\]
 \end{prop}
 This is the refinement of Proposition \ref{prop_S_box_triv} that we sought.
 
 \subsection{Proof of Theorem \ref{thm_gen_PPW}}

 We will now input this bound for $S_4(\taubf)$ with the choice $\Kbf = 2\Hbf/\Pbf$ into Proposition \ref{prop_TS4}, always with the specification that $r>M$ and $\Qbf$ is chosen as in (\ref{Q_choice}).
For each $i=1,\ldots, k$ we choose $P_i$  such that 
\[\frac{1}{2} H_iq_i^{- \frac{1}{2(r-M)}} \leq P_i < H_iq_i^{-\frac{1}{2(r-M)}}.\]
With this choice, (\ref{Kmin'}) is satisfied; we also have $P_i \geq 1$ as long as $H_i > q_i^{\frac{1}{2(r-M)}}$.
In order to satisfy the conditions $H_i P_i < q_i$ of Lemma \ref{lemma_A}, we must also restrict to ranges with $H_i < q_i^{\frac{1}{2} + \frac{1}{4(r-M)}}$. With these choices, we apply Proposition \ref{prop_S_box_nontriv} in Proposition \ref{prop_TS4} to obtain
\[  T(\Fbf;\Nbf,\Hbf) \ll  \| \Hbf\|^{-\frac{1}{2r}} \| \Pbf \|^{1 - \frac{1}{2r}} \left \{ \Qbf^\ga K_{\max}^{2rk+\ep}K_{\min}^{-M} \| \qbf \|^{1/2}  \right\}^{\frac{1}{2r}}.
\] 
We recall that
\[ \Qbf^{\frac{\ga}{2r}} \ll (\Hbf/\Pbf)^{\frac{\ga}{2r}} \ll  (\qbf^\ga)^{\frac{1}{4r(r-M)}} .\]
Thus we may conclude
\begin{eqnarray*}
  T(\Fbf;\Nbf,\Hbf)& \ll&  \| \Hbf\|^{-\frac{1}{2r}} \| \Pbf \|^{1 - \frac{1}{2r}} (\qbf^\ga)^{\frac{1}{4r(r-M)}} q_{\max}^{\frac{2rk}{4r(r-M)}}  q_{\min}^{-\frac{M}{4r(r-M)}} \| \qbf \|^{\frac{1}{4r}+\ep}  \\
  & \ll&   \| \Hbf\|^{1-\frac{1}{r}} \| \qbf \|^{ - \frac{1}{2(r-M)}(1- \frac{1}{2r})} (\qbf^\ga)^{\frac{1}{4r(r-M)}} q_{\max}^{\frac{2rk}{4r(r-M)}}  q_{\min}^{-\frac{M}{4r(r-M)}} \| \qbf \|^{\frac{1}{4r}+\ep}  \\
& \ll &  \| \Hbf\|^{1-\frac{1}{r}} \| \qbf \|^{ \frac{-r-M+1}{4r(r-M)}+\ep} (\qbf^\ga)^{\frac{1}{4r(r-M)}} q_{\max}^{\frac{2rk}{4r(r-M)}}  q_{\min}^{-\frac{M}{4r(r-M)}}  ,
  \end{eqnarray*}
which proves Theorem \ref{thm_gen_PPW}. (Here we note that certainly $r> R(d+1)$ implies $r>M$.)

In the case where $q_i=q$ for all $i$, we have $q_{\max} = q_{\min} =q$, $\|\qbf \| = q^k$, and $\qbf^\ga = q^M$ where $M$ is the weight of the system $\Fbf$, so that this simplifies to 
\[ T(\Fbf;\Nbf,\Hbf) \ll
 \| \Hbf\|^{1-\frac{1}{r}} q^{\frac{k(r+1-M)}{4r(r-M)} + \ep}.\]
 
\subsection{Optimal choice of $r$}\label{sec_thm_Vin_proof}
Suppose that for each $i$, $H_i = q^{1/4 + \kappa_i}$. Set $\sig = \sum_{i=1}^k \kappa_i.$ Then Corollary \ref{thm_gen_PPW_corq} provides an upper bound of size $\| H \| q^{-\del + \ep}$ where 
\[ \del = \frac{4\sig(r-M) - k}{4r(r-M)}.\]
\xtra{
We write
\[ \del = \frac{\sig}{r} - \frac{k}{4r(r-M)} = \frac{c}{x} - \frac{1}{ax(x-b)} = f(x),\]
where $a=4/k, b=M, c=\sig$,
so that the derivative is
\[ f'(x) = -cx^{-2} + (ax^2-abx)^{-2}(2ax-ab). \]
Now $f'(x)$ has a root precisely at the solutions 
\[ cax^2 + (-2abc-2)x +ab^2c+b=0.\]
We compute the solutions to this are 
\[ x=\frac{abc+1\pm \sqrt{abc+1}}{ac} = b + \frac{1\pm \sqrt{abc+1}}{ac}  .\]
Since $f(x)$ blows down to $-\infty$ at $x=b$ (assuming $c$ small relative to $a$), we see that the root we are interested in takes $+$ instead of minus. Thus we choose 
\[ r = b + \frac{1 + \sqrt{abc+1}}{ac}  .\]
}
As a function of $r$, this attains a maximum at the real value 
\[ r_0 =M + \frac{k\left(1 + \sqrt{\frac{4M \sig}{k} +1}\right)}{4\sig}.\]
Choosing $r = r_0 + \theta$ where $-1/2 \leq \theta < 1/2$, we see that $\del$ is approximately of size 
\[ \del \approx \frac{4\sig^2}{k \left( 1+ \sqrt{ \frac{4M\sig }{k} + 1}\right)^2}.\]
For fixed $k,d$ as $\sig= \sum \kappa_i \maps 0$ this behaves like 
\[
\del \approx \frac{ \sig^2}{k},
\]
which we note is nicely dependent only on the dimension $k$ and not on other parameters of the system $\Fbf$.

 \section{Technical lemmas}\label{sec_tech}
 \subsection{Proof of Lemma \ref{lemma_A}}\label{sec_A_lemma}
 We now return to the proof of Lemma \ref{lemma_A}. It is clear from the definition of $\Acal(\mbf)$ that it vanishes unless each $m_i$ satisfies $|m_i| \leq 2q_i$. Next we note that $\Acal(\mbf)$ is a non-negative integer, so trivially $S_1 \leq S_2$. Thus we turn to bounding $S_2$, for which we note that 
 \begin{multline*}
 \sum_\mbf \Acal(\mbf)^2 
 	 = 
 \sum_\mbf \# \{ \pbf, \pbf', \abf ,\abf' : 0 \leq a_i < p_i, 0 \leq a_i' < p_i' \; \text{and} \; p_i,p_i' \in \Pcal_i:  \\
 m_i  \leq  \frac{N_i - a_iq_i}{p_i} < m_i + \frac{H_i}{P_i} ,  m_i  \leq  \frac{N_i - a_i'q_i}{p_i'} < m_i + \frac{H_i}{P_i}  , i=1,\ldots, k \} .
 \end{multline*}
For a fixed $\mbf$, in order for a quadruple $\pbf,\pbf',\abf,\abf'$ to belong to this set we must have both $(N_i - a_iq_i)/p_i$ and $(N_i - a_i'q_i)/p_i'$ belong to the interval $[m_i, m_i + H_i/P_i)$ (for all $i$), so that we require
\[  \left| \frac{N_i - a_iq_i}{p_i} - \frac{N_i - a_i'q_i}{p_i'} \right| \leq \frac{H_i}{P_i}, \quad \text{for each $i=1,\ldots, k$}.\]
If these conditions are satisfied then there will be $O(\|\Hbf \| \, \|\Pbf \|^{-1})$ corresponding values $\mbf$ for which this can occur.
We may thus deduce that 
\begin{eqnarray}
  \sum_\mbf \Acal(\mbf)^2 
  &\ll& \|\Hbf \| \, \|\Pbf \|^{-1} \# \{ \pbf, \pbf', \abf, \abf'  : 0 \leq   \left| \frac{N_i - a_iq_i}{p_i} - \frac{N_i - a_i'q_i}{p_i'} \right| \leq \frac{H_i}{P_i} \} \nonumber \\
  &\ll& \|\Hbf \| \, \|\Pbf \|^{-1} \sum_{\pbf, \pbf' \in \Pscr} \Mcal(\pbf,\pbf'), \label{AMPP}
  \end{eqnarray}
  where we set
  \begin{multline*}
   \Mcal(\pbf,\pbf') :=\# \{ \mathbf{0} \leq \abf < \pbf, \mathbf{0} \leq \abf'  <\pbf' : \\
    0 \leq   \left| \frac{N_i - a_iq_i}{p_i} - \frac{N_i - a_i'q_i}{p_i'} \right| \leq \frac{H_i}{P_i},  \text{for each $i=1,\ldots, k$} \} .
    \end{multline*}
   We now define for any primes $p_i,p_i' \in \Pcal_i$ the quantity  
  \[ M_i(p_i,p_i') = \# \{0 \leq a_i<p_i, 0 \leq a_i'< p_i':   0 \leq   \left| \frac{N_i - a_iq_i}{p_i} - \frac{N_i - a_i'q_i}{p_i'} \right| \leq \frac{H_i}{P_i}\} .\] 
  We note that for each pair of tuples $\pbf, \pbf'$,
\[  \Mcal (\pbf,\pbf') = \prod_{i=1}^k M_i(p_i, p_i').\]
Thus 
\begin{eqnarray}
 \sum_{\pbf, \pbf' \in \Pscr} \Mcal (\pbf,\pbf')  
 	& = &  \sum_{\pbf, \pbf' \in \Pscr} \left( \prod_{i=1}^k M_i(p_i, p_i') \right)  \nonumber \\
 &=& \prod_{i=1}^k \left( \sum_{p_i, p_i' \in \Pcal_i} M_i(p_i,p_i') \right) \nonumber \\
	 &=&  \prod_{i=1}^k \left( \sum_{p_i = p_i' \in \Pcal_i} M_i(p_i,p_i')  + \sum_{p_i  \neq p_i' \in \Pcal_i} M_i(p_i,p_i')  \label{MPP} \right) .
	 \end{eqnarray}
We now recall from the proof of Lemma 2.2 in \cite{HBP14a} that in the one-dimensional case it is already known that for each $i=1,\ldots, k$,
\begin{eqnarray*}
\sum_{p_i \in \Pcal_i} M_i(p_i,p_i) & \ll&  P_i^2\\
\sum_{p_i \neq p_i' \in \Pcal_i} M_i(p_i,p_i') & \ll & P_i^2 ,
\end{eqnarray*}
with the latter bound holding under the condition $P_iH_i< q_i$.
Applying this in (\ref{MPP}), we see that 
\[  \sum_{\pbf, \pbf' \in \Pscr} \Mcal (\pbf,\pbf')  \ll \| \Pbf \|^{2},\]
so that in total (\ref{AMPP}) shows that
\[ \sum_\mbf \Acal(\mbf)^2 \ll \| \Hbf \| \|\Pbf \|,\]
as desired.

\subsection{Proof of Lemma \ref{lemma_BoIw}}\label{sec_BoIw}
Recall that in Lemma \ref{lemma_BoIw} we consider the sum 
\beq\label{I_sum_0}
 \sum_{\nbf \in I } a(\nbf),
 \eeq
for arbitrary complex numbers $a(\nbf)$ indexed by $\nbf \in \Z^k$ lying in an arbitrary fixed product of sub-intervals $I \subseteq (\Abf, \Abf+\Bbf]$. We will denote $I = (\Cbf, \Cbf + \Dbf]$, with $(C_i,C_i+D_i] \subseteq (A_i,A_i+B_i]$ for each $i=1,\ldots, k$. We note that if any $D_i=0$, then the sum (\ref{I_sum_0}) is vacuous; thus we may assume all $D_i>0$. Next note that if $B_i<1$ then there is at most one value $n_i$ considered in the $i$-th coordinate of the sum $\sum_{\nbf \in (\Abf,\Abf+\Bbf]}a(\nbf)$, and we could regard the sum as living in a lower dimensional setting and proceed with the proof in a lower dimension. Thus we may assume $B_i \geq 1$ for all $i = 1,\ldots, k$.

We will prove Lemma \ref{lemma_BoIw} with a simple adaptation of Bombieri and Iwaniec's original argument \cite{BoIw86}. For each $i$, let $\psi_i(x)$ denote a $C^\infty$ compactly supported non-negative function that vanishes for $x \leq \lfloor C_i \rfloor$ and $x \geq \lfloor C_i+D_i \rfloor+1$ and is identically 1 for $\lfloor C_i \rfloor +1 \leq x \leq \lfloor C_i+D_i \rfloor$; clearly we may also choose this so that $|\psi_i| \leq 1$ and $\psi_i$ has uniformly bounded derivatives $|\psi_i^{(N)}| \leq 1$ for all $N \geq 1$.
Let $\Psi(\xbf) = \psi_1(x_1) \cdots \psi_k(x_k)$, so that 
\[ \sum_{\nbf \in (\Cbf, \Cbf+\Dbf]} a(\nbf) = \sum_{\nbf \in (\Abf, \Abf + \Bbf]} a(\nbf)\Psi(\nbf).\]
After expressing $\Psi(\xbf)$ in terms of its inverse Fourier transform (see (\ref{Psi_L1})), we have
\begin{eqnarray*}
\sum_{\nbf \in (\Cbf, \Cbf + \Dbf]} a(\nbf)
	&=& \sum_{\nbf \in (\Abf, \Abf + \Bbf]} a(\nbf) \int_{\R^k} \hat{\Psi}(\boldsymbol{\theta}) e(\nbf \cdot \boldsymbol{\theta}) d \boldsymbol{\theta} \nonumber \\
	& = &   \int_{\R^k} \hat{\Psi}(\boldsymbol{\theta}) \sum_{\nbf \in (\Abf, \Abf + \Bbf]} a(\nbf) e(\nbf \cdot \boldsymbol{\theta}) d \boldsymbol{\theta} .
	\end{eqnarray*}
Thus
\beq\label{anpsin}
\left| \sum_{\nbf \in (\Cbf, \Cbf + \Dbf]} a(\nbf) \right|
	 \leq  \sup_{\theta \in \R^k} \left|\sum_{\nbf \in (\Abf, \Abf + \Bbf]} a(\nbf) e(\nbf \cdot \boldsymbol{\theta}) \right| \|\hat{\Psi} \|_{L^1(\R^k)}, 
\eeq
	where
	\[   \|\hat{\Psi} \|_{L^1(\R^k)} = \int_{\R^k} |\hat{\Psi}(\boldsymbol{\theta})| d \boldsymbol{\theta}.\]
We now note that 
\[ \hat{\Psi}(\boldsymbol{\theta}) = \int_{\R^k} \Psi(\xbf)e(-\xbf \cdot \boldsymbol{\theta}) d\xbf  = \prod_{i=1}^k \left(  \int_{A_i-1}^{A_i+B_i+1} \psi_i(x_i)e(-x_i\theta_i) dx_i \right) = \prod_{i=1}^kJ_i(\theta_i),\]
say, where we define
\[ J_i(\theta) = \int_{A_i-1}^{A_i+B_i+1} \psi_i(x)e(-x\theta) dx .\]
Each of these satisfies 
\beq\label{int_triv}
|J_i(\theta)|  \ll \min\{ B_i+2, |\theta|^{-N}\}, \quad \text{for any $N \geq 1$}.
\eeq
The first option is the trivial bound; the second option follows from integration by parts. Precisely, for a fixed $\theta \neq 0$, by writing
\[ e^{-2\pi i \theta x} = \frac{1}{(-2\pi i \theta)} \frac{d}{dx} e^{-2\pi i \theta x},\]
we may integrate by parts repeatedly to see that 
\[\int_{A_i-1}^{A_i+B_i+1} \psi_i(x)e(-x\theta) dx  =  \frac{1}{(2\pi i \theta)^N}\int_{A_i-1}^{A_i+B_i+1} \psi_i^{(N)}(x)e(-x\theta) dx \]
for any $N \geq 1$; the boundary terms vanish due to the compact support of $\psi_i$. Now we note that $\psi_i^{(N)}(x)$ is uniformly bounded by assumption, and moreover it vanishes unless $x$ belongs to either of two intervals of length $1$. 
This gives the desired result (\ref{int_triv}).

We temporarily set $L_i = B_i+2$ and apply (\ref{int_triv}) to observe that for each $i$:
 \begin{eqnarray*}
 \int_{-\infty}^\infty |J_i (\theta) | d\theta
 	& \ll & \int_{|\theta| \leq L_i^{-1}} L_i d \theta + \int_{L_i^{-1} \leq |\theta| \leq 2L_i} |\theta|^{-1} d \theta + \int_{|\theta| \geq 2L_i} |\theta|^{-2} d\theta \nonumber \\
	& \ll & 1 + \log (L_i) + L_i^{-1}  \ll \log (B_i+2).
	\end{eqnarray*}
	
Finally, we see that 
\[  \int_{\R^k} |\hat{\Psi}(\boldsymbol{\theta})| d \boldsymbol{\theta}
	= \prod_{i=1}^k \left( \int_{-\infty}^\infty |J_i (\theta_i) | d\theta_i \right)  \ll \prod_{i=1}^k \log (B_i+2),\]
confirming that 
\beq\label{Psi_L1}
\hat{\Psi} \in L^1 (\R^k).
\eeq
 Thus the use of the Fourier inversion formula is justified, and we also see in (\ref{anpsin}) that 
\[ \left| \sum_{\nbf \in (\Cbf, \Cbf + \Dbf]} a(\nbf)\right|  \ll ( \prod_{i=1}^k \log (B_i+2) )  \sup_{\theta \in \R^k} \left|\sum_{\nbf \in (\Abf, \Abf + \Bbf]} a(\nbf) e(\nbf \cdot \boldsymbol{\theta}) \right|, \]
as claimed.
	
\subsection{Proof of Lemma \ref{lemma_part_sum_gen}}\label{sec_part_sum_proof}
We will proceed by iterated partial summation applied to
\beq\label{ab_sum}
\sum_{\boldsymbol{0} < \nbf \leq \Nbf} a(\nbf) b(\nbf) .
\eeq
We first apply partial summation with respect to $n_1$ in (\ref{ab_sum}). We set $J=\{1\}$ and $I=\{2,\ldots, k\}$  so that
\begin{eqnarray*}
 \sum_{\nbf\leq \Nbf} a(\nbf) b(\nbf) &=& 
	\sum_{\nbf_{(I)} \leq \Nbf_{(I)} }\left( \sum_{0<n_1 \leq N_1} a(n_1, \nbf_{(I)}) b(n_1 ,\nbf_{(I)}) \right)\\
	&=&\sum_{\nbf_{(I)} \leq \Nbf_{(I)}} \left\{ b(N_1,\nbf_{(I)}) \left( \sum_{0<n_1 \leq N_1} a(n_1, \nbf_{(I)}) \right)  \right.\\
		&& \qquad \qquad -\, \left. \int_0^{N_1} \left( \sum_{0<n_1 \leq t_1} a(n_1, \nbf_{(I)}) \right) \frac{d}{dt_1}  b(t_1, \nbf_{(I)}) dt_1 \right\}.
	\end{eqnarray*} 
We may then apply partial summation with respect to $n_2$, and so on, iteratively for each $n_i$ with $i \leq k$.
One obtains a representation of (\ref{ab_sum}) as a sum of $2^k$ terms, each corresponding to a subset $J\subseteq \{ 1, \ldots, k\}$ (and its corresponding complement $I$). For each partition $J\union I$ of $\{1,\ldots, k\}$ with $|J| = v$, the resulting term is of the shape
\[(-1)^{v} \idotsint_{\bstack{(0,N_j]}{j \in J}} A_{(I),(J)}(\Nbf_{(I)},\tbf_{(J)}) \frac{ \partial^{v}}{\partial \tbf_{(J)}} b(\Nbf_{(I)},\tbf_{(J)})  d\tbf_{(J)}.\]
 Here if $J = \{j_1,\ldots, j_v\}$ we let $ \frac{ \partial^{v}}{\partial \tbf_{(J)}}= \frac{ \partial^{v}}{\partial t_{j_1} \cdots t_{j_v}}$ and $d\tbf_{(J)} = dt_{j_1} \cdots dt_{j_v}$.
As a result of the assumed bounds (\ref{b_assp}) on the partial derivatives of $b(\xbf)$, we may conclude that 
\[ \left| \sum_{\nbf\leq \Nbf} a(\nbf) b(\nbf) \right| \leq \sum_{J \subseteq \{1,\ldots, k\}} \left( \prod_{j \in J} B_j \right)
 \idotsint_{\bstack{(0,N_j]}{j \in J}} \left| A_{(^cJ),(J)}(\Nbf_{(^cJ)},\tbf_{(J)}) \right| d\tbf_{(J)},\]
 which is the statement of the lemma.

\subsection*{Acknowledgments}
The author is supported in part by NSF DMS-1402121, and also thanks the Hausdorff Center for Mathematics for a very pleasant working environment during a portion of this research. The author thanks D. R. Heath-Brown and T. D. Wooley for informative discussions on the work of this paper, D. Schindler for a close reading of an earlier version of the manuscript, and the anonymous referees for helpful comments on the exposition.

\bibliographystyle{amsplain}
\bibliography{NoThBibliography}

\end{document}

%% file: format.tex
\newtheorem{letterthm}{Theorem}

\newtheorem{thm}{Theorem}[section]
\newtheorem*{thm*}{Theorem}
\newtheorem{prop}[thm]{Proposition}
\newtheorem{lemma}[thm]{Lemma}
\newtheorem{cor}{Corollary}[thm]

\newcommand{\ep}{\varepsilon}

\newcommand{\con}{\equiv}

\newcommand{\ndiv}{\nmid}
\newcommand{\modd}[1]{\; ( \mathrm{mod} \; #1)}
\newcommand{\bstack}[2]{#1 \atop #2}

\newcommand{\maps}{\rightarrow}
\newcommand{\intersect}{\cap}

\newcommand{\union}{\cup}

\newcommand{\al}{\alpha}
\newcommand{\be}{\beta}

\newcommand{\ga}{\gamma}
\newcommand{\del}{\delta}
\newcommand{\Del}{\Delta}
\newcommand{\ka}{\kappa}

\newcommand{\Sig}{\Sigma}
\newcommand{\sig}{\sigma}

\newcommand{\Acal}{\mathcal{A}}
\newcommand{\Pcal}{\mathcal{P}}

\newcommand{\Bcal}{\mathcal{B}}
\newcommand{\Dcal}{\mathcal{D}}

\newcommand{\Mcal}{\mathcal{M}}

\newcommand{\Scal}{\mathcal{S}}

\newcommand{\Abf}{\mathbf{A}}
\newcommand{\Bbf}{\mathbf{B}}
\newcommand{\Lbf}{\mathbf{L}}
\newcommand{\Hbf}{\mathbf{H}}
\newcommand{\Mbf}{\mathbf{M}}
\newcommand{\Qbf}{\mathbf{Q}}
\newcommand{\Rbf}{\mathbf{R}}

\newcommand{\Ubf}{\mathbf{U}}

\newcommand{\Dbf}{\mathbf{D}}

\newcommand{\Nbf}{\mathbf{N}}
\newcommand{\Kbf}{\mathbf{K}}

\newcommand{\Xbf}{\mathbf{X}}

\newcommand{\abf}{{\bf a}}

\newcommand{\mbf}{{\bf m}}
\newcommand{\nbf}{{\bf n}}

\newcommand{\sbf}{{\bf s}}
\newcommand{\tbf}{{\bf t}}

\newcommand{\xbf}{{\bf x}}
\newcommand{\ybf}{{\bf y}}
\newcommand{\zbf}{{\bf z}}

\newcommand{\F}{\mathbb{F}}

\newcommand{\N}{\mathbb{N}}

\newcommand{\R}{\mathbb{R}}

\newcommand{\Z}{\mathbb{Z}}

\newcommand{\Lscr}{\mathscr{L}}
\newcommand{\Pscr}{\mathscr{P}}

\newcommand{\beq}{\begin{equation}}
\newcommand{\eeq}{\end{equation}}